# Complex surfaces of general type: some recent progress

Ingrid C. Bauer[1], Fabrizio Catanese[2], and Roberto Pignatelli[3]⋆

[1] Mathematisches Institut, Lehrstuhl Mathematik VIII, Universitätstraße 30, D-95447 Bayreuth, Germany `Ingrid.Bauer@uni-bayreuth.de`
[2] Mathematisches Institut, Lehrstuhl Mathematik VIII, Universitätstraße 30, D-95447 Bayreuth, Germany `Fabrizio.Catanese@uni-bayreuth.de`
[3] Dipartimento di Matematica, Università di Trento, via Sommarive 14, I-38050 Povo (TN), Italy `Roberto.Pignatelli@unitn.it`

**Introduction**

In this article we shall give an overview of some recent developments in the theory of complex algebraic surfaces of general type.

After the *rough* or *Enriques - Kodaira* classification of complex (algebraic) surfaces, dividing compact complex surfaces in four classes according to their Kodaira dimension $-\infty$, 0, 1, 2, the first three classes nowadays are quite well understood, whereas even after decades of very active research on the third class, the class of surfaces *of general type*, there is still a huge number of very hard questions left open. Of course, we made some selection, which is based on the research interest of the authors and we claim in no way completeness of our treatment. We apologize in advance for omitting various very interesting and active areas in the theory of surfaces of general type as well as for not being able to mention all the results and developments which are important in the topics we have chosen.

Complex surfaces of general type come up with certain (topological, birational) invariants, topological as for example the topological Euler number $e$ and the self intersection number of the canonical divisor $K^2$ of a minimal surface, which are linked by several (in-) equalities. In the first chapter we will summarize the classically known inequalities, which force surfaces of general type in a certain region of the plane having $K^2$ and $e$ as coordinates, and we shall briefly comment on the so-called *geography* problem, whether, given numerical invariants lying in the admissible range, i.e., fulfilling the required inequalities, does there exist a surfaces having these invariants. We shall however more broadly consider the three classical invariants $K^2, p_g, q$, which determine the other invariants $\chi := 1 - q + p_g, e = 12\chi - K^2$.

⋆ The present work was performed in the realm of the SCHWERPUNKT "Globale Methoden in der komplexen Geometrie", and was also supported by a VIGONI-DAAD Program. A first draft of this article took origin from the lectures by the second author at the G.A.C. Luminy Meeting , october 2005: thanks to the organizers!

An important new inequality, which Severi tried without success to establish, and which has been attacked for many years with partial results by several authors, asserts that a surface of maximal Albanese dimension satisfies the inequality $K^2 \geq 4\chi$. We will report on Pardini's surprisingly simple proof of this so-called Severi's conjecture (cf. [Par05]).

The study of the pluricanonical maps is an essential technique in the classification of surface of general type. The main results concerning the $m$-canonical maps with $m \geq 3$ go back to an earlier period and we refer to [Cat87b] for a report on them.

We will report in the second chapter on recent developments concerning the bicanonical map; we would like to mention Ciliberto's survey (cf. [Cil97]) on this topic for the state of art ten years ago. Here instead, we combine a discussion of this topic with the closely intertwined problem of classification of surfaces with low values of the numerical invariants.

In the third chapter we report on surfaces of general type with geometric genus $p_g$ equal to four, a class of surfaces whose investigation was started by Federigo Enriques (cf. chapter VIII of his book 'Le superficie algebriche', [Enr49]).

By Gieseker's theorem we know that for fixed $K^2$ and $\chi$ there exists a quasi projective coarse moduli space $\mathcal{M}_{K^2,\chi}$ for the birational equivalence classes of surfaces of general type. It is a very challenging problem to understand the geometry of these moduli spaces even for low values of the invariants. The case $p_g = 4$ is studied via the behaviour of the canonical map. While it is still possible to divide the moduli space into various locally closed strata according to the behaviour of the canonical map, it is very hard to decide how these strata patch together.

Using certain presentations of Gorenstein rings of codimension 4 introduced by M. Reid and D. Dicks, which arrange the defining equations as Pfaffians of certain matrices with many symmetries in such a way that these equations behave well under deformation, it is possible to exhibit explicit deformations, which allow to "connect" certain irreducible components of the moduli space.

Inspired by a construction of A. Beauville of a surface with $K^2 = 8$, $p_g = q = 0$, the second author defined Beauville surfaces as surfaces which are rigid and which admit an unramified covering which is isomorphic to a product of curves of genus at least 2. In this case the moduli space of surfaces orientedly homeomorphic to a given surface consists either of a unique real point, or of a pair of complex conjugate points corresponding to complex conjugate surfaces.

These surfaces, and the more general surfaces isogenous to a product, not only provide cheap counterexamples to the Friedman - Morgan speculation (which will be treated more extensively in the sixth section of this article), but provide also a wide class of surfaces quite manageable in order to test conjectures, and offer also counterexamples to various problems. The ease with which one can handle these surfaces is based on the fact that these surfaces are determined by "discrete" combinatorial data.

Beauville surfaces, their relations to group theory and to Grothendieck's theory of 'Dessins d'enfants' will be discussed in the fourth chapter.

It is a very difficult and very intriguing problem to decide whether two algebraic surfaces, which are not deformation equivalent, are in fact diffeomorphic.

The theory of Lefschetz fibrations provides an algebraic tool to prove that two surfaces are diffeomorphic. By a theorem of Kas (which holds also in the symplectic context) two Lefschetz fibrations are diffeomorphic if and only if their corresponding factorizations of the identity in the mapping class group are equivalent under the equivalence relation generated by Hurwitz moves and by simultaneous conjugation. We outline the theory, which was used with success in [CW04] in chapter five, which we end with a brief report on the status of two very old conjectures by Chisini concerning cuspidal curves and algebraic braids.

As already mentioned before, one of the fundamental problems in the theory of surfaces of general type is to understand their moduli spaces, in particular the connected components which parametrize the deformation equivalence classes of minimal surfaces of general type. By a classical result of Ehresmann, two deformation equivalent algebraic varieties are diffeomorphic. The other direction, i.e., whether two diffeomorphic minimal surfaces of general type are indeed in the same connected component of the moduli space, was an open problem since the eighties. We discuss in the last chapter the various counterexamples to the Friedman - Morgan speculation, who expected a positive answer to the question (unlike the second author, cf.[Kat83]).

Moreover, we briefly report on another equivalence relation introduced by the second author, the so - called quasi étale-deformation (Q.E.D.) equivalence relation, i.e., the equivalence relation generated by birational equivalence, by quasi étale morphisms and by deformation equivalence. For curves and surfaces of special type two varieties are Q.E.D. equivalent if and only if they have the same Kodaira dimension, whereas there are infinitely many surfaces of general type, which are pairwise not Q.E.D. equivalent.

# 1 Old and new inequalities

## 1.1 Invariants of surfaces

Let $X$ be a compact complex manifold and let $\Omega_X^n$ be its canonical bundle, i.e., the line bundle of holomorphic $n$−forms (usually denoted by $\omega_X$, since it is a dualizing sheaf in the sense of Serre duality). A corresponding canonical divisor is usually denoted by $K_X$.

To $X$ one associates its *canonical ring*

$$R(X) := \oplus_{m \geq 0} H^0(\omega_X^{\otimes m}).$$

The trascendency degree over $\mathbb{C}$ of this ring leads to

- the *Kodaira dimension* $\kappa(X) := tr(R(X)) - 1$,

if $R(X) \neq \mathbb{C}$, otherwise $\kappa(X) := -\infty$. The Kodaira dimension is invariant under deformation (by Siu's theorem [Siu02], generalizing Iitaka's theorem for surfaces) and can assume the values $-\infty, 0, \ldots, n = \dim X$.

**Definition 1.** $X$ *is said to be* of general type *if the Kodaira dimension is maximal,* $\kappa(X) = \dim X$.

We are interested in the case of *surfaces*, i.e., of manifolds of dimension 2, of general type.

The three principal invariants under deformations for the study of these surfaces are

- the self intersection of the canonical class $K_S^2$ of a minimal model,
- the geometric genus $p_g := h^0(\omega_X)$ and
- the irregularity $q := h^1(\mathcal{O}_S) = h^0(\Omega_S^1)$.

The equality $h^1(\mathcal{O}_S) = h^0(\Omega_S^1)$ follows by Hodge theory since every algebraic surface is projective.

The invariants we have introduced, with the exception of $K_S^2$, are not only deformation invariants but also birational invariants.

**Definition 2.** *A smooth surface $S$ is called* minimal *(or* a minimal model*) iff it does not contain any exceptional curve $E$ of the first kind (i.e. $E \cong \mathbb{P}^1$, $E^2 = -1$).*

Every surface can be obtained by a minimal one (its "minimal model") after a finite sequence of blowing ups of smooth points; this model is moreover unique if $\kappa(S) \geq 0$ (see III.4.4, III.4.5 and III.4.6 of [BHPV04]). Thus, every birational class of surfaces of general type contains exactly one minimal surface, and one classifies surfaces of general type by studying their minimal models. To each minimal surface of general type we will associate its numerical

- type $(K_S^2, p_g, q)$,

a triple of integers given by the three invariants introduced above.

In fact these determine all other classical invariants, as

- the Euler-Poincaré characteristic of the trivial sheaf $\chi(\mathcal{O}_S) = 1 - q + p_g$;
- the topological Euler characteristic $e(S) = c_2(S) = 12\chi(\mathcal{O}_S) - K_S^2$;
- the plurigenera $P_m(S) := h^0(\omega_X^{\otimes m}) = \chi(\mathcal{O}_S) + \binom{m}{2}K_S^2$.

The epression for $c_2$ is a classical theorem of M. Noether, and the expression for the plurigenera follows by Riemann-Roch and by Mumford's vanishing theorem.

By the theorems on pluricanonical maps (cf.[Bom73]), minimal surfaces $S$ of general type with fixed invariants are birationally mapped to normal surfaces $X$ in a fixed projective space of dimension $P_5(S) - 1$. $X$ is uniquely determined, is called the *canonical model* of $S$, and is obtained contracting to points all the (-2)-curves of $S$ (curves $E \cong \mathbb{P}^1$, with $E^2 = -2$).

Let us recall Gieseker's theorem

**Theorem 1 (Gieseker [Gie77]).** *There exists a quasi-projective coarse moduli scheme for canonical models of surfaces of general type $S$ with fixed $K_S^2$ and $c_2(S)$.*

In particular, we can consider the subscheme $\mathcal{M}_{K_S^2, p_g, q}$ corresponding to minimal surfaces of general type of type $(K_S^2, p_g, q)$. By the above theorem, it is a quasi projective scheme, in particular, it has finitely many irreducible components.

It is a dream ever since to completely describe $\mathcal{M}_{K_S^2, p_g, q}$ for as many types as possible.

### 1.2 Classical inequalities and geography

Obviously the first question is: for which values of $(K_S^2, p_g, q)$ is $\mathcal{M}_{K_S^2, p_g, q}$ non empty?

For example, it is clear that $p_g(S)$ and $q(S)$ are always nonnegative, since they are dimensions of vector spaces.

In fact much more is known. In the following table we collect the well known classical inequalities holding among the invariants of minimal surfaces of general type:

$$\begin{array}{lll} & K_S^2 \geq 1 & \chi \geq 1 \\ (N) & K_S^2 \geq 2p_g - 4 & \text{or the weaker } K_S^2 \geq 2\chi(\mathcal{O}_S) - 6 \\ (D) & \text{if } q > 0, \ K_S^2 \geq 2p_g & \text{or the weaker if } q > 0, \ K_S^2 \geq 2\chi(\mathcal{O}_S) \\ (MY) & K_S^2 \leq 9\chi & \end{array}$$

We have labeled by (N)= Noether, (D) = Debarre, (MY) = Miyaoka-Yau the rows, corresponding to the names of the inequalities ([Deb82], [Deb83], [Miy77], [Yau78], see also [BHPV04], chap. 7).

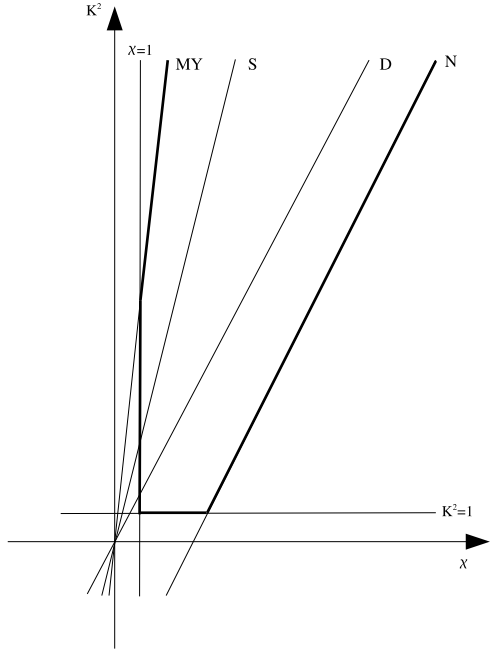

**Fig. 1.** The geography of minimal surfaces of general type

In figure 1 we have drawn the limit lines (i.e., where equality holds) of the various inequalities in the $(\chi, K_S^2)$ - plane.

The above listed inequalities show that the pair of invariants $\chi, K_S^2$ of a minimal surface of general type gives a point with integral coordinates in the convex region limited by the "bold" piecewise linear curve. Moreover, if $q > 0$ this point cannot be at the "right" of the line $D$.

We drew one more line in our picture, labeled by $(S)$. This is the Severi line $K^2 = 4\chi$, i.e., the equality case of the Severi inequality $K^2 \geq 4\chi \Leftrightarrow K^2 \geq \frac{1}{2}e$, which will be discussed in detail at the end of this section.

### 1.3 Surfaces fibred over a curve

An important method for the study of surfaces of general type is to consider relatively minimal fibrations of surfaces over curves $f : S \to B$.

**Definition 3.** *A fibration $f : S \to B$ is a surjective morphism with connected fibres. We are interested in the case of fibrations* of surfaces to curves, *meaning that in this paper $S$ and $B$ will always be smooth compact complex manifolds of respective dimensions* 2 *and* 1.

*The fibration is said to be* relatively minimal *if $f$ does not contract any rational curve of self intersection $-1$ to a point.*

One denotes

- by $b$ the genus of the base curve $B$;
- by $g$ the genus of a general fibre.

To avoid confusion, let us point out that a fibration is called rational or irrational according to the genus $b$ of the base being $0$ or $> 0$. On the other hand, the genus of the fibration is the genus $g$ of the fibre. For example, if we say $f$ is a genus 2 rational fibration, we intend that $g = 2$ and $b = 0$.

The classical way of saying: a genus $b$ pencil of curves of genus $g$ is however still the most convenient way to describe a fibration.

To a relatively minimal fibration $f$ one associates

- its relative canonical bundle $\omega_{S|B} := \omega_S \otimes f^*(\omega_B^\vee)$ and
- the sheaves ($\forall n \geq 0$) $V_n := f_*(\omega_{S|B}^{\otimes n})$.

The sheaves $V_n$ are vector bundles (i.e., locally free sheaves) with very nice properties.

**Theorem 2 (Fujita [Fuj78a], [Fuj78b]).** *The vector bundles $V_n$ are semi-positive, i.e., every locally free quotient of it has nonnegative degree.*

To be more precise, $V_1$ is a direct sum of an ample vector bundle with $q(S) - b$ copies of the trivial bundle and with some undecomposable stable degree 0 vector bundle without global sections. Zucconi [Zuc97] proved moreover that if one of those stable bundles has rank 1, then it is a torsion line bundle.

For $n \geq 2$ we have:

**Theorem 3 (Esnault-Viehweg [EV90]).** $\forall n \geq 2$ *the vector bundle $V_n$ is ample unless $f$ has constant moduli, which means that all the smooth fibres are isomorphic.*

Since $R^1 f_* \omega_{S|B} = \mathcal{O}_B$ by relative duality, and $R^1 f_* \omega_{S|B}^{\otimes n} = 0$ $\forall n \geq 2$ by the assumption of relative minimality, one can compute the Euler characteristic of $V_n$ by Riemann-Roch, and consequently its degree.

We introduce the following invariants of the fibration $f$:

- the self intersection of the relative canonical divisor

$$K_f^2 := \omega_{S|B} \cdot \omega_{S|B} = K_S^2 - 8(g-1)(b-1),$$

- the Euler characteristic of the relative canonical divisor

$$\chi_f = \chi(\omega_{S|B}) = \chi(\mathcal{O}_S) - (g-1)(b-1),$$

- its slope $\lambda(f) := K_f^2/\chi_f$.

The slope is clearly defined only for $\chi_f \neq 0$, or equivalently (as we will see soon) if the fibration is not a holomorphic bundle.

The above mentioned computation gives

$$\deg V_n = \chi_f + \frac{n(n-1)}{2} K_f^2$$

and since by Fujita's theorem these numbers are nonnegative this gives the two inequalities $K_f^2 \geq 0$ and $\chi_f \geq 0$ respectively known as Arakelov's inequality (cf. [Ara71]) and Beauville's inequality (cf. [Bea82]).

In fact, we have the following list of inequalities

(A)   $K_f^2 \geq 0$, i.e., $K_S^2 \geq 8(g-1)(b-1)$,
(B)   $\chi_f \geq 0$, i.e., $\chi(\mathcal{O}_S) \geq (g-1)(b-1)$,
(ZS) $c_2(S) \geq 4(b-1)(g-1)$,
(NN) $q \leq b + g$,
(X)   $4 - \frac{4}{g} \leq \lambda(f) \leq 12$.

Here the meaning of the labeling is the following: (A) = Arakelov's inequality, (B) = Beauville' inequality, (X) = Xiao's inequality (also known as slope inequality), (NN) = no name's inequality, (ZS) = Zeuthen-Segre. A proof of those inequalities can be found in [Bea82] with the exception of the slope inequality, proved in [Xia87] (see also [CH88] in the semistable case).

The equality cases of the first 4 inequalities are well described:

- if equality holds in (A), $f$ has constant moduli;
- equality holds in (B) $\Leftrightarrow$ $f$ has constant moduli and is smooth;
- for $g \geq 2$, equality holds in (ZS) $\Leftrightarrow$ $f$ is smooth;
- $q = b + g \Leftrightarrow f$ is birationally equivalent to the projection of a product $B \times F$ to the first factor.

In particular, we see that the slope is defined whenever the fibration is not a holomorphic bundle, since the denominator $\chi_f$ vanishes iff equality holds in Beauville's inequality.

An important consequence is the following

**Theorem 4 (Beauville).** *If $X$ is a minimal surface of general type, then $p_g \geq 2q - 4$. Moreover, if $p_g = 2q - 4$, then $S$ is a product of a curve of genus 2 with a curve of genus $q - 2$.*

Note for later use (see next section) the following

**Corollary 1.** *If $p_g = q$ (i.e., if $\chi(\mathcal{O}_S) = 1$), then $p_g = q \leq 4$. Moreover, minimal surfaces of general type with $p_g = q = 4$ are exactly the products of two genus 2 curves.*

*Proof of theorem 4.* The standard wedge product on 1−forms induces a natural map

$$\wedge : \Lambda^2 H^0(\Omega_S^1) \to H^0(\Omega_S^2)$$

Recall that $q = \dim H^0(\Omega^1_S)$, $p_g = \dim H^0(\Omega^2_S)$. Let us assume $p_g \leq 2q - 4$.

By a dimension count, if $p_g \leq 2q - 4$, the projective linear subspace of $\mathbb{P}(\Lambda^2 H^0(\Omega^1_S))$ corresponding to the kernel of the above map must intersect the Plücker embedding of the Grasmannian $G_2(H^0(\Omega^1_S))$ (which has dimension $2q - 4$), and therefore there are two linearly independent 1−forms $\omega_1$ and $\omega_2$ such that the following holomorphic two form is identically zero: $\omega_1 \wedge \omega_2 \equiv 0$.

By the theorem of Castelnuovo-De Franchis there is a fibration $f : S \to B$ with base of genus $b \geq 2$, and two holomorphic 1−forms $\alpha_1, \alpha_2 \in H^0(\Omega^1_B)$ such that $f^*\alpha_i = \omega_i$. Since $S$ is of general type, also $g \geq 2$.

Then

$$\chi_f \geq 0 \Rightarrow \chi(\mathcal{O}_S) \geq (b-1)(g-1) = (b-2)(g-2) + b + g - 3 \geq q - 3.$$

So we have $1 - q + p_g \geq q - 3 \Leftrightarrow p_g \geq 2q - 4$.

If $p_g = 2q - 4$, all inequalities are equalities and then, since $q = b + g$ and $(b-2)(g-2) = 0$, $S$ is a product of two curves of genus at least 2, and one of the two must have genus exactly 2.

□

## 1.4 Severi's inequality

We recall that the Albanese variety $Alb(X)$ of a compact Kähler manifold $X$ is the cokernel of the natural map

$$\int : H_1(X, \mathbb{Z}) \to H^0(\Omega^1(X))^\vee$$

defined by integrating 1−forms on 1-cycles.

The Albanese morphism

$$\alpha : X \to Alb(X)$$

is defined (up to translations in $Alb(X)$) by fixing a point $p_0 \in X$, and by associating to each point $p \in X$ the class in $Alb(X)$ of $\int_{p_0}^p$, where the integral is taken along any path between $p_0$ and $p$.

Recall that, if $X$ is projective (as any surface of general type), $Alb(X)$ is an abelian variety (of dimension $q$).

The Albanese morphism is a powerful tool for studying *irregular* surfaces ($q > 0$) and in particular:

**Definition 4.** *A variety $X$ is called of* maximal Albanese dimension *if the image of the Albanese morphism has the same dimension as $X$.*

This is the general case for surfaces, since otherwise the Albanese morphism is a fibration onto a smooth curve of genus $q$. We see then that for surfaces maximal Albanese dimension is equivalent to the non existence of a genus $q$ pencil.

We can now state the theorem known as *Severi's inequality*

**Theorem 5 (Pardini [Par05]).** *If $S$ is a smooth complex minimal surface of maximal Albanese dimension, then $K_S^2 \geq 4\chi$.*

This theorem was proved only very recently by R. Pardini, but it has a long story, which we briefly sketch in the following.

### Severi's conjecture

The inequality takes its name from F. Severi, since he was the first to claim the result in the 30's [Sev32].

His proof turned out to be wrong, as was pointed out in [Cat83], since it was based on the assertion that a surface with irregularity $q$ either contains an irrational genus $q$ fibration, or the sections of $H^0(\Omega_S^1)$ have no common zero. Counterexamples were given in [Cat84], where there were constructed bidouble covers $S \to X$ of any algebraic surface with, among other properties, $q(S) = q(X)$. If $X$ has no irrational pencils, since the Albanese map of $S$ factors through the cover, then also $S$ has no irrational pencils. But any ramification point of the cover is a base point for $H^0(\Omega_S^1)$.

Therefore Severi's inequality was posed in [Cat83] as *Severi's conjecture*, a conjecture on surfaces of general type, since for surfaces with $\kappa(S) \leq 1$ it is a straightforward consequence of the Enriques-Kodaira classification. It had also been posed as a conjecture by M. Reid (conj. 4 in [Rei79]) who proved the weaker $K_S^2 \geq 3\chi$.

### Proofs in special cases

In the 80's, Xiao's work on surfaces fibred over a curve was mainly motivated by Severi's conjecture. In [Xia87] he proved the slope inequality and Severi's conjecture for surfaces having an irrational pencil.

In the 90's Konno [Kon96] proved the conjecture in the special case of *even* surfaces, i.e., surfaces whose canonical class is $2-$ divisible in the Picard group.

Finally, at the end of the 90's, Manetti [Man03] could prove the inequality for surfaces of general type whose canonical bundle is ample.

### Manetti's proof

Manetti considers the tautological line bundle $L$ of the $\mathbb{P}^1-$bundle $\pi : \mathbb{P}(\Omega_S^1) \to S$; standard computations give

$$3(K_S^2 - 4\chi) = L^2 \cdot (L + \pi^* K_S).$$

Then, using the fact that $\Omega_S^1$ is generically globally generated, he can write the right hand side of the above equation as $2K_S E + (L + \pi^* K_S)C$ for an effective $1-$cycle $C$ in $\mathbb{P}(\Omega_S^1)$, and where $E$ is the maximal effective divisor in $S$ such that $h^0(\Omega_S^1(-E)) = h^0(\Omega_S^1)$. Thus the problem is reduced to the nonnegativity of the term $(L+\pi^*K_S)C$. This is obvious if $\Omega^1(K_S)$ is nef, but in general it requires a very detailed and complicated analysis of the $1$−cycle $C$.

In fact, Pardini's proof not only does not require the ampleness of the canonical divisor, but is much easier than Manetti's.

We should however mention that Manetti's argument leads to a very detailed description of the equality case, showing that a surface of general type of maximal Albanese dimension, lying on the Severi line ($K^2 = 4\chi$), and having ample canonical class, has irregularity $q = 2$ and is a double cover of a principally polarized Abelian surface, branched on a divisor $D$ algebraically equivalent to $2\Theta$.

Up to now there is no similar description of the limit case without the assumption that $K$ be ample.

### Pardini's proof

Pardini's idea is to construct a sequence of genus $g_d$ fibrations $f_d : Y_d \to \mathbb{P}^1$ such that
$$\lim_{d \to \infty} g_d = +\infty \text{ and } \lim_{d \to \infty} \lambda(f_d) = K_S^2/\chi(\mathcal{O}_S).$$

Then, taking the limit of the left-hand side of the slope inequality, one gets the desired inequality $K_S^2/\chi(\mathcal{O}_S) \geq 4$.

To construct these fibrations, she considers the Cartesian diagram

$$\begin{array}{ccc} S' & \xrightarrow{p} & S \\ \alpha' \downarrow & & \downarrow \alpha \\ Alb(S) & \xrightarrow{\cdot d} & Alb(S) \end{array},$$

where $d : Alb(S) \to Alb(S)$ is multiplication by $d$.

One observes that $S'$ is connected since we have a surjection $\pi_1(S) \to \pi_1(Alb(S)) = H_1(S, \mathbb{Z})$.

Clearly, $K_{S'}^2 = d^{2q}K_S^2$, $\chi(\mathcal{O}_{S'}) = d^{2q}\chi(\mathcal{O}_S)$.

Let $L$ be a very ample divisor on $Alb(S)$ and set $H := \alpha^*L$, $H' := \alpha'^*L$. Then $p^*H \sim_{num} d^2 H'$, whence $H'^2 = d^{2q-4}H^2$ and $K_{S'}H' = d^{2q-2}K_S H$.

Let now $D_1, D_2 \in |H'|$ be two general curves and define $C_1 := D_1 + D_2 \in |H_1 + H_2|$. Moreover, choose $C_2 \in |2H'|$ sufficiently general such that $C_1$ and $C_2$ intersect transversally. $C_1$ and $C_2$ define a rational pencil $f_d : Y_d \to \mathbb{P}^1$, where $Y_d$ is the blow up of $S'$ at $C_1 \cap C_2$. The singular fibre induced by $C_1$ guarantees that $f_d$ is not a holomorphic bundle, whence the slope $\lambda(f_d)$ is well defined.

For the invariants of $Y_d$ we get
$$K_{Y_d}^2 = K_{S'}^2 - 4H'^2 = d^{2q}K_S^2 - 4d^{2q-4}H^2$$

$$\chi(Y_d) = \chi(S') = d^{2q}\chi(S)$$

$$g_d = 1 + K_{S'}H' + 2H'^2 = 1 + d^{2q-2}K_S H + 2d^{2q-4}H^2$$

and therefore $\lim_d g_d = +\infty$ as requested.

Moreover, $K^2_{f_d} = K^2_{Y_d} + 8(g_d - 1)$ and $\chi(f_d) = \chi(Y_d) + (g_d - 1)$ and we see that both invariants are polynomials in $d$ of degree $2q$, whose leading terms are respectively $K^2_S$ and $\chi(\mathcal{O}_S)$. In particular,

$$\lim_{d\to\infty} \lambda(f_d) = \lim_{d\to\infty} K^2_f/\chi_f = K^2_S/\chi(\mathcal{O}_S)$$

## 2 Surfaces with $\chi = 1$ and the bicanonical map

### 2.1 The bicanonical map

The behaviour of the $m - th$ canonical map of $S$ (i.e., the rational map associated to $|mK_S|$) is an essential tool in the theory of surfaces of general type.

As we mentioned in the introduction, the cases where $m \geq 3$ are solved since long (cf. the survey [Cat87b]).

The canonical map ($m = 1$) was first studied by Beauville in [Bea79], but there remain still many unresolved questions.

The case $m = 2$ was particularly studied in the last years, and we have the impression that we are very close to a complete understanding.

In order to fix the starting point, we summarize the results of several authors ([Fra91],[Rei88],[Cat81],[CC91],[CC93],[Xia85a]) in the following

**Theorem 6.** *Let $S$ be a minimal surface of general type. Then*

- *the bicanonical map is generically finite unless $p_g = 0$ and $K^2 = 1$;*
- *if $K^2_S \geq 5$ or $p_g \geq 1$, the bicanonical map is a morphism.*

Note that, if $p_g = 0$ and $K^2 = 1$, then $P_2 = 2$ and the bicanonical map is a rational ($b = 0$) fibration. In all known examples this is a genus 4 fibration, although at the moment it is only proven that its genus is 3 or 4 (see [CP05]).

These surfaces are usually called *numerical Godeaux surfaces*. Numerical Godeaux surfaces with torsion (in the Picard group) of cardinality at least 3 are classified in [Rei78], a family with torsion $\mathbb{Z}/2$ was constructed in [Bar84]. Up to last year only sporadic examples of surfaces with trivial torsion were known, but recently Schreyer [Sch05] has announced the construction of a family of the expected dimension ($= 8$) using a new approach based on homological algebra.

The above theorem says that in all other cases the bicanonical map maps $S$ to a surface, and it is a morphism (except for finitely many families).

In the last years many people studied the degree of this map, in particular, trying to classify the surfaces such that the bicanonical map is not birational.

**The standard case**

It is well known that the bicanonical map of a smooth curve of general type (i.e., of genus at least 2) fails to be birational if and only if the curve has genus 2.

This exception induces a *"standard exception"* to the birationality of the bicanonical map in dimension 2.

**Definition 5.** *A surface $S$ of general type* presents the standard case *if there exists a dominant rational map onto a curve $f : S \dashrightarrow B$ whose general fibre is irreducible of genus* 2.

In fact, if $S$ presents the standard case, then the restriction of the bicanonical map of $S$ to a general fibre factors through the bicanonical map of the fibre itself and therefore cannot be birational.

The subschemes of the moduli space corresponding to surfaces of presenting the standard case are not empty for infinitely many moduli spaces, and Persson [Per81] constructed many interesting surfaces considering double covers of ruled surfaces branched on relative sextics, thereby filling a big region of the convex region represented in figure 1.

Bombieri ([Bom73]) showed that the standard case gives almost all exceptions to the birationality of the bicanonical map. More precisely, combining his results with those of Reider ([Rei88]) we know now that a minimal surface of general type with $K^2 \geq 10$ either presents the standard case, or its bicanonical map is birational. In particular, the exceptions to the birationality of the bicanonical map not presenting the standard case belong to finitely many families and many authors are trying since then to classify them.

**du Val's double planes**

In the same paper [Bom73] Bombieri constructed a surface of type $(K^2, p_g, q) = (9, 6, 0)$ not presenting the standard case. His example can be easily described as a hypersurface $F_{14}$ of degree 14 in the weighted projective space $\mathbb{P}(1, 1, 2, 7)$, and from this description it follows rightaway that it bicanonical map is a double cover of the weighted projective space $\mathbb{P}(1, 1, 2)$ (isomorphic to a quadric cone in $\mathbb{P}^3$). This example is in fact a special case of a more general "geometric" situation studied first by du Val.

Let $S$ be a minimal regular surface with $p_g \geq 2$, such that the general canonical curve is irreducible, smooth and hyperelliptic. Since the restriction of the bicanonical map $\varphi_{2K}$ to a canonical curve factors through the canonical map of the curve itself, $\varphi_{2K}$ cannot be birational.

Du Val [Duv52] gave a list of such surfaces obtained as double covers of rational surfaces. A generalization (see [CML00], [Bor03]) leads to the following:

**Definition 6 (du Val's double planes).** *A smooth surface $S$ is* a du Val double plane *if it is birational to*

$\mathcal{D}$) a double cover of $\mathbb{P}^2$ branched over a smooth curve of degree $8$;

$\mathcal{D}_n$) a double cover of $\mathbb{P}^2$ branched over the union of a curve of degree $10 + n$ with $n$ distinct lines through a point $p$, such that the essential singularities of the branch curve are the following:
- $p$ is a singular point of multiplicity $2n + 2$,
- there is a singular point of type $[5, 5]$ on each line,
- possibly there are some quadruple points and some points of type $[3, 3]$;

$\mathcal{B}$) a double cover of the Hirzebruch surface $\mathbb{F}_2$ whose branch curve can be decomposed as $C_0 + G'$, $G' \in |7C_0 + 14\Gamma|$ (where $|\Gamma|$ is the ruling of $\mathbb{F}_2$ and $C_0$ is the section with self intersection $-2$), whose only essential singularities are $[3, 3]$ points that are tangent to a fibre.

Recall that a singular point of type $[d, d]$ is a singular point of multiplicity $d$ having a further singular point of multiplicity $d$ infinitely near to the first one. In other words, if we blow-up the singular point, the strict transform of the curve has one more singular point of multiplicity $d$ lying on the exceptional divisor.

*Remark 1.* In the definition of the du Val's double planes of type $\mathcal{D}_n$ and $\mathcal{B}$ we only care about the essential singularities of the branch curve (as usual in the theory of double covers) since adding a simple singularity to the branch curve does not affect the properties of the resulting surface we are interested in.

On the contrary, in the definition of the double planes of type $\mathcal{D}$, we assume the branch curve to be smooth. In fact, if we take a double cover of $\mathbb{P}^2$ branched over a curve of degree $8$ with a double point, the pull back of the pencil of lines through this point to the surface defines a pencil of curves of genus $2$ (through a singular point of the surface), so the resulting surface presents the standard case.

Note that this example shows that one can degenerate surfaces presenting a nonstandard case to surfaces presenting the standard case (just take a family of smooth plane curves of degree $8$ degenerating to a singular one and consider the corresponding family of double covers).

Borrelli proved that this list is "complete" in the following sense

**Theorem 7 (Borrelli [Bor03]).** *If $S$ is a minimal surface of general type, not presenting the standard case, whose bicanonical map factors through a degree $2$ rational map onto a rational or ruled surface: then $S$ is the smooth minimal model of a du Val double plane. In particular, either $q = 0$ or $p_g = q = 1$.*

### The "classification"

The standard case and the du Val's double planes do not give all possible surfaces of general type with nonbirational bicanonical map, but the remaining exceptions are really few.

What is known about these is summarized in the following

**Theorem 8.** *Let $S$ be a smooth minimal surface of general type whose bicanonical map $\varphi_{2K}$ is not birational. Then one of the following cases occur:*

i)  *$S$ presents the standard case;*
ii) *$S$ is the smooth minimal model of a du Val double plane;*
iii) *$S$ is a surface of type $(1,1,0)$ (for these automatically $\deg \varphi_{2K} = 4$ and $S$ is a complete intersection of two sextics in $\mathbb{P}(1,2,2,3,3)$);*
iv) *$S$ is of type $(2,1,0)$ with Picard group having torsion $\mathbb{Z}/2\mathbb{Z}$ (for these automatically $\deg \varphi_{2K} = 4$ and its double cover corresponding to the torsion class is a complete intersection of two quartics in $\mathbb{P}(1,1,1,2,2)$);*
v)  *$\varphi_{2K}$ is $2:1$ onto a K3 surface and $p_g = 1$, $q = 0$, $2 \leq K^2 \leq 8$;*
vi) *$S$ is of type $(6,3,3)$ or of type $(4,2,2)$ (for these both cases automatically $\varphi_{2K}$ has degree 2 and we have a nonstandard case)*
vii) *$S$ has $p_g = q \leq 1$.*

All these cases with the exception of $p_g = q \leq 1$ are now rather clear.

The history of this theorem is rather complicated and combines the efforts of several authors. We try to reconstruct its more important steps here, giving some more details on each class.

One of the first results in this direction is due to Xiao Gang [Xia90], giving, in the nonstandard case, and under the assumption that the degree $d$ of the bicanonical map is at least 3, a list of the possible values of $d$ and of the possible places in the Enriques classification of the bicanonical image $\Sigma$.

In 1997, Ciliberto, Francia and Mendes Lopes [CFML97] gave a complete classification of the case $p_g \geq 4$, essentially confirming du Val's list.

Then, Ciliberto and Mendes Lopes, with contributions of the second author and Borrelli, worked in the next years to extend the classification to $p_g \geq 2$ (see [CCML98], [CML00], [CML02a], [CML02b], [Bor02]). The case $p_g = 1$ and $q = 0$, giving cases $iii)$, $iv)$ and $v)$ is classified in [Bor03].

In fact, cases $iii)$ and $iv)$ resulted already from the analysis of [Xia90] where it is proven that, if $\deg \varphi_{2K} \geq 3$, then either $S$ is of type $(1,1,0)$ or of type $(2,1,0)$, or with $p_g = q \leq 2$.

The description given in $iii)$ and $iv)$ of the first two cases comes from the papers [Cat79] and [CD89], where all surfaces of respective types $(1,1,0)$ and $(2,1,0)$ are classified. In particular, it is shown that all surfaces of type $(1,1,0)$ are as in $iii)$.

*Remark 2.* Surfaces of type $(2,1,0)$ without torsion in homology, also sometimes called Catanese-Debarre surfaces, offer the following interesting phenomenon: there is an irreducible component of the moduli space such that

1) for the general surface the bicanonical map is birational, while there are subvarieties for which the bicanonical map can respectively be

2) of degree 2 onto a K3- quartic surface,
3) of degree 2 onto a rational quartic surface,
4) of degree 4 onto a smooth quadric surface.

The surfaces in $v)$ are usually called Todorov surfaces, since they were introduced in [Tod81]. The subspaces of the moduli spaces corresponding to them is described in [Mor88].

Finally, the largely open case $vii)$ is very strongly related with the problem, of independent interest, of the classification of surfaces of general type with $p_g = q$.

We shall describe in the next subsection what is currently known on these surfaces, showing in particular that we have a very precise description of the two cases in $vi)$.

### 2.2 Surfaces with $p_g = q$

These are the surfaces corresponding to the "vertical" piece of the bold line in figure 1. In particular, $1 \leq K^2 \leq 9$.

#### Surfaces with $p_g = q \geq 4$.

This case is clear, by corollary 1 of Beauville's theorem 4. If $p_g = q \geq 4$, then $S$ is a product of two genus 2 curves and $p_g = q = 4$. We recall that then $K^2 = 8$ and clearly the bicanonical map has degree 4, and we have a standard case.

#### Surfaces with $p_g = q = 3$.

These surfaces have been first studied in [CCML98], and a complete classification has been recently achieved independently by Pirola [Pir02] and Hacon-Pardini [HP02].

The result is the following

**Theorem 9.** *A minimal surface of general type with $p_g = q = 3$ has $K^2 = 6$ or $K^2 = 8$ and, more precisely,*

- *if $K^2 = 6$, $S$ is the symmetric square of a genus $3$ curve;*
- *otherwise $S = C_2 \times C_3/\tau$, where $C_g$ denotes a curve of genus $g$ and $\tau$ is an involution of product type acting on $C_2$ as an elliptic involution (i.e., with elliptic quotient), and on $C_3$ as a fixed point free involution.*

*In particular, the moduli space of minimal surfaces of general type with $p_g = q = 3$ is the disjoint union of $\mathcal{M}_{6,3,3}$ and $\mathcal{M}_{8,3,3}$, which are both irreducible of respective dimension $6$ and $5$.*

We sketch the idea of the proof.

By Debarre's inequality (in the "stronger" form: $q > 0 \Rightarrow K_S^2 \geq 2p_g$), $p_g = q = 3$ implies $K^2 \geq 6$.

As in the proof of Beauville's theorem, consider now the map

$$\wedge : \Lambda^2(H^0(\Omega_S^1)) \to H^0(\Omega_S^2).$$

Since $p_g = q = 3$, it is a linear map between two three dimensional spaces. If this is not an isomorphism, then (since every vector in $\Lambda^2\mathbb{C}^3$ is decomposable) there are two nontrivial $1-$forms $\omega_1$ and $\omega_2$ with $\omega_1 \wedge \omega_2 \equiv 0$. This yields then (by Castelnuovo-De Franchis) a pencil $f : S \to B$ with $b, g \geq 2$.

Then, by Beauville's inequality, $b = g = 2$ and the fibration is a holomorphic bundle (this forces $K_S^2 = 8$). Therefore $f$ is induced by a map $\pi_1(B) \to Aut(F)$ (where $F$ is a smooth fibre), whose kernel induces an unramified cover $\varphi : C \to B$ Galois with group $G$.

Since $q = 3$, the quotient of $F$ by the group $G$ has genus 1. By Hurwitz's formula one easily sees that, if $\phi : F \to F/G$ is branched in one point, then $|G| \leq 4$, hence $G$ is Abelian, contradicting that $\phi$ is ramified. Again Hurwitz's formula shows that $\phi$ is branched in 2 points and $G \cong \mathbb{Z}/2$.

Otherwise $\wedge$ is an isomorphism, and therefore $S$ does not have any pencil $f : S \to B$ with $b \geq 2$. In particular $\alpha(S)$ is a surface, a divisor $\Theta$ in $Alb(S)$. Pirola noticed that $\Theta$ must be ample, else it would have an elliptic fibration and therefore an irrational pencil with base of genus $b \geq 2$.

This implies, by Lefschetz's hyperplane theorem, that the induced map $H^1(\Omega_{Alb(S)}^1) \to H^1(\Omega_S^1)$ is injective: since, for any class $\eta \in H^1(\Omega_{Alb(S)}^1)$ $\int c_1(\Theta) \wedge \eta \wedge \bar{\eta} > 0$.

In particular, $h^1(\Omega_S^1) \geq 9$ and this (since by Hodge theory $12\chi - K^2 = c_2 = 2p_g - 4q + 2 + h^1(\Omega^1)$) implies $K_S^2 \leq 7$.

The case $K_S^2 = 6$ was already settled in [CCML98], where it is first shown that the degree of the scheme of base points is $K^2 - 6$, and then that in the case $K_S^2 = 6$ $\alpha$ is an embedding. More precisely it is shown that its image is a theta divisor in a principally polarized abelian threefold and therefore $S$ is the symmetric square of a genus 3 curve.

What remains to prove is $\mathcal{M}_{7,3,3} = \emptyset$, and this is done in [Pir02] by a careful study of the paracanonical system.

The fact that the bicanonical map has degree 2 is an easy consequence of the adjunction formula by which $K_S$ is the pull back of $\Theta$, and of the fact that the sections of $\mathcal{O}_A(2\Theta)$ are invariant, as well as $\Theta$, for the symmetry of $A$ sending $x \to -x$.

**Surfaces with $p_g = q = 2$.**

This case is still far from being classified. Ciliberto and Mendes Lopes [CML02a] classified all surfaces with $p_g = q = 2$ and non-birational bicanonical map (not presenting the standard case). Their result, corresponding to the subcase $(4, 2, 2)$ of case vi) of theorem 8, is the following

**Theorem 10.** *If $S$ is a minimal surface of general type with $p_g = q = 2$ and non-birational bicanonical map not presenting the standard case, then $S$ is a double cover of a principally polarized abelian surface $(A, \Theta)$, with $\Theta$* *irreducible. The double cover $S \to A$ is branched along a divisor $B \in |2\Theta|$, having at most double points. In particular $K_S^2 = 4$.*

Note that, again by Debarre's inequality, $p_g = q = 2 \Rightarrow K^2 \geq 4$, so Ciliberto and Mendes Lopes' surfaces belong to the limit case. Their theorem solves completely the problem of the non birationality of the bicanonical map in this case, but of course a complete classification of minimal surfaces of general type with $p_g = q = 2$ would be interesting by itself.

Results in this direction have been recently obtained by F. Zucconi; to explain them we need to give the following definition.

**Definition 7.** *A surface $S$ is said to be isogenous to a (higher) product if $S$ admits an unramified finite covering which is biholomorphic to a product of two curves of respective genera at least 2.*

We have already seen surfaces isogenous to a product in our analysis of surfaces with $p_g = q$, namely all surfaces in $\mathcal{M}_{8,4,4}$ and all surfaces in $\mathcal{M}_{8,3,3}$.

Zucconi's theorem is the following

**Theorem 11 (2.9 in [Zuc03]).** *There are two classes of minimal surfaces of general type with $p_g = q = 2$ whose Albanese image is a surface and having an irrational pencil, and they are both isogenous to a higher product.*

*More precisely, either they have a double cover which is a product of two genus 2 curves or they are a quotient of the product of two genus 3 curves by an action of $\mathbb{Z}/2\mathbb{Z}$.*

In both cases Zucconi describes precisely the group action as a diagonal action induced by actions on the two curves. The interested reader will find all details in Zucconi's paper.

Zucconi managed also to remove the hypothesis on the Albanese map, by use of a special class of surfaces isogenous to a higher product, the generalized hyperelliptic surfaces introduced in [Cat00].

**Definition 8.** *Let $C_1$ and $C_2$ be two smooth curves, $G$ a finite group with two injections respectively in $Aut(C_1)$ and $Aut(C_2)$. Then the quotient surface $S = C_1 \times C_2/G$ by the diagonal action is said to be a generalized hyperelliptic surface if*

- *the projection $C_1 \to C_1/G$ is unramified;*
- *$C_2/G$ is rational.*

Then Zucconi proved

**Theorem 12.** *If $S$ has $p_g = q = 2$, and the image of the Albanese map is a curve, then $S$ is a generalized hyperelliptic surface.*

What remains to be classified is the class of surfaces with $p_g = q = 2$ having no irrational pencils.

Chen and Hacon, in a preprint, constructed an example of surfaces with $p_g = q = 2$, $K^2 = 5$ and Albanese morphism of degree 3.

**Surfaces with $p_g = q = 1$.**

In this case Debarre's inequality gives only $K^2 \geq 2$.

The Albanese morphism is a map onto an elliptic curve, in particular, all these surfaces have a fibration with base of genus $b = 1$. We summarize in the following statement what is known about these surfaces.

**Theorem 13.**

- $\mathcal{M}_{2,1,1}$ is unirational (by this we mean: irreducible and unirational) of dimension 7. The Albanese map of all these surfaces is a genus 2 fibration.
- $\mathcal{M}_{3,1,1}$ has 4 connected components, all unirational of dimension 5. The Albanese map is a genus 3 fibration for the surfaces in one of those components, and a genus 2 fibration in all other cases.
- $\mathcal{M}_{4,1,1}$, $\mathcal{M}_{5,1,1}$ and $\mathcal{M}_{8,1,1}$ are non empty.

Actually much more can be said, and we try to be more precise in the following.

First the most mysterious cases. It remains unsettled the existence of surfaces of general type with $p_g = q = 1$ and $K^2 = 6, 7, 9$. In a recent preprint by Rito appears the construction of a surface with $p_g = q = 1$ and $K^2 = 6$ as a double cover of a Kummer surface modifying slightly (i.e., adding a singular point to the branch curve) Todorov's construction of a surface with $p_g = 1, q = 0$ and $K^2 = 8$ in [Tod81]. Its construction makes use of the computer program MAGMA (to find a branch curve with the right singularities).

Second, the "partially understood" cases. Examples of surfaces with $p_g = q = 1$ and $K^2 = 4, 5$ were constructed by the second author as bidouble covers in [Cat99]. In both cases the Albanese map turns out to be a genus 2 fibration, so they present the standard case. The case $K^2 = 8$ was studied by Polizzi [Pol06], who considered the cases of surfaces having bicanonical map of degree 2. He could prove that all these surfaces are isogenous to a product and that they form three components of the moduli space, one of dimension 5 and two of dimension 4. All these surfaces do not contain any genus 2 pencil and they are in fact du Val double planes.

Finally, the cases $K^2 = 2, 3$ are completely classified.

The first to be settled was $K^2 = 2$, done by the second author in [Cat81], representing all those surfaces as double covers of the symmetric square of their Albanese curve.

The case $K^2 = 3$ was first studied in [CC91] where it was shown, among other things, that the Albanese map could be either a genus 2 or a genus 3 fibration. The case $g = 3$ was then classified in [CC93], showing that it gives a unirational family of dimension 5.

Note that, if there is surface with $p_g = q = 1$, $K^2 \leq 3$ and nonbirational bicanonical map not presenting the standard case, it must belong to this family.

The question whether such a surface exists is still open. Recently Polizzi [Pol05] has shown that a general surface in this component has birational bicanonical map, but this is not true for all of them, since Xiao[Xia85b] has found a subfamily of dimension 1 having a genus 2 pencil.

The classification of the case $K^2 = 3$ was completed in [CP05] classifying all those surfaces having Albanese fibres of genus 2.

The main tool for this classification is a new method for studying fibrations $f : S \to B$ of genus 2, and fibrations of genus 3 with general fibre non hyperelliptic, basically giving *generators* and *relations* of their relative canonical algebra $\mathcal{R}(f) = \oplus V_n$, seen as a sheaf of algebras over $B$.

Let us recall the vector bundles $V_n$ introduced in the previous section as $V_n = f_* \omega_{S|B}^{\otimes n}$. Roughly speaking then, $\mathcal{R}(f)$ is a bundle whose fibres are the canonical rings of the fibres of $f$.

We state here only the theorem for genus 2 fibrations, since it is the one used in order to complete this classification.

**Theorem 14.** *A genus* 2 *fibration* $f : S \to B$ *is determined by the following* 5 *data*

- *the base curve* $B$;
- *the rank* 2 *vector bundle* $V_1 := f_* \omega_{S|B}$ *over* $B$;
- *an effective divisor* $\tau$ *on* $B$;
- *a class* $\xi \in Ext^1_{\mathcal{O}_B}(S^2(V_1), \mathcal{O}_\tau)/(Aut_{\mathcal{O}_B}(\mathcal{O}_\tau)$ *yielding* $V_2$;
- *letting* $\mathcal{A}$ *be the subring of the relative canonical algebra generated by* $V_2$, $V_3^+$ *the* (+1) *eigenbundle for the hyperelliptic involution on the fibres, and defining* $\tilde{\mathcal{A}}_6 := \mathcal{H}om((V_3^+)^2, \mathcal{A}_6)$ *(where* $\mathcal{A}_6$ *is the image in* $\mathcal{A}$ *of* $S^3(V_2)$*), the last datum is an element* $w \in \mathbb{P}(H^0(\tilde{\mathcal{A}}_6))$.

*Moreover,* $\deg V_1 = \chi(\mathcal{O}_S) - (b-1)$, $\deg \tau = K_S^2 - 2\chi(\mathcal{O}_S) - 10(b-1)$.

We want to explain here the geometry behind this theorem, which at a first glance can appear slightly technical.

The vector bundles $V_n$ yield the degree $n$ part of the canonical ring of each fibre. So each of these vector bundles induces a rational map, the *relative n-canonical map*, from $S$ to the corresponding projective bundle $\mathbb{P}(V_n)$, mapping each fibre via its $n$-canonical map.

The multiplication map of degree 1 forms give a morphism of sheaves $S^2(V_1) \to V_2$ which fits into an exact sequence

$$0 \to S^2(V_1) \to V_2 \to \mathcal{O}_\tau \to 0$$

for an effective divisor $\tau$ on $B$ supported on the image of the "bad" fibres (those which are not 2-connected, i.e., the fibres that can be decomposed as $A + B$ with $A, B$ effective divisors such that $A \cdot B = 1$).

$\xi$ is the class of this extension. Therefore $\xi$ yields $V_2$ and determines the relative bicanonical map. Since the bicanonical map of a genus 2 curve is a double cover of a conic (branched in 6 points), this map has degree 2 onto a conic subbundle $\mathcal{C}$ of the $\mathbb{P}^2$-bundle $\mathbb{P}(V_2)$. We have $\mathcal{C} = Proj(\mathcal{A})$.

In [CP05] it is proven that the relative bicanonical map is a morphism contracting at most some rational curves with self-intersection $(-2)$, which implies that the branch curve has no "essential" singularities.

In fact the $5^{th}$ datum $w$ determines the branch curve $div_{\mathcal{C}}(w)$. In fact, sections of $S^3(V_2)$, or of a twist of it, are "equations" of divisors in $\mathbb{P}(V_2)$ which cut a cubic curve on each fibre. Taking the quotient by the subsheaf corresponding to the equations vanishing on the conic bundle $\mathcal{C}$, one gets an equation for a divisor on the conic bundle, which cuts 6 points (intersection of a cubic and a conic) on a general fibre, and gives our branch curve.

Let us come back to the case $p_g = q = 1$, $K^2 = 3$ and $g = 2$. One needs to construct a suitable genus 2 fibration over an elliptic curve, with, in the sense of theorem 14, $\deg V_1 = \deg \tau = 1$. This is done in [CP05], by studying vector bundles on elliptic curves, and three different families are found.

Let us finally mention that $\mathbb{P}(V_1)$ is the symmetric square of $B$. In fact, our double cover is birational to a double cover of it. The behaviour of this double cover was described in [CC91], characterising these surfaces as double covers of the symmetric product of an elliptic curve with branch locus belonging to a certain algebraic system with prescribed singularities.

This new method shows then, rather surprisingly, that this algebraic system is not connected.

### Surfaces with $p_g = q = 0$.

The class of surfaces with $p_g = q = 0$ is one of the most complicated and intriguing classes of surfaces of general type. By the standard inequalities we have: $1 \leq K_S^2 \leq 9$.

We have already mentioned the case $K^2 = 1$, of the numerical Godeaux surfaces, the only case for which the bicanonical map is not finite, so let us restrict to $K_S^2 \geq 2$.

These surfaces are very far from being classified. From the point of view of the bicanonical system, this case was object of an intensive analysis by Mendes Lopes and Pardini in the last years.

What it is known on the degree of the bicanonical map can be summarized in the following

**Theorem 15 ([MLP05],[MLP02]).** *Let $S$ be a surface with $p_g = q = 0$. Then*

- *if $K^2 = 9 \Rightarrow \deg \varphi_{2K} = 1$,*
- *if $K^2 = 7, 8 \Rightarrow \deg \varphi_{2K} = 1$ or $2$,*
- *if $K^2 = 5, 6 \Rightarrow \deg \varphi_{2K} = 1, 2$ or $4$,*
- *if $K^2 = 3, 4 \Rightarrow \deg \varphi_{2K} \leq 5$ and if moreover $\varphi_{2K}$ is a morphism, then $\deg \varphi_{2K} = 1, 2$ or $4$,*

- if $K^2 = 2$ then obviously, since the image of the bicanonical map is $\mathbb{P}^2$, then the bicanonical map is non birational, and obviously we have: $\deg \varphi_{2K} \leq 8$, equality holding if and only if $\varphi_{2K}$ is a morphism.

We would like to mention that, by Reider's theorem, the bicanonical map is a morphism as soon as $K^2 \geq 5$. In fact the bicanonical map of all known examples of surfaces with $p_g = q = 0$ and $K^2 \geq 2$ is a morphism. So one could suspect that the bicanonical map is always a morphism whenever $K^2 \geq 2$.

Langer ([Lan00]) has proven that the bicanonical map of a minimal surface of general type with $p_g = q = 0$ and $K_S^2 = 4$ has no fixed part.

Mendes Lopes and Pardini gave also a description of some of these surfaces having non birational bicanonical map, in particular for $K^2 \geq 6$. We summarize here some of their results

**Theorem 16 ([MLP03],[MLP01],[MLP04a],[MLP04b]).** *Let $S$ be a minimal surface of general type with $p_g = q = 0$ whose bicanonical map is not birational. Then the image of the bicanonical map is a rational surface unless $K^2 = 3$, $\deg \varphi_{2K} = 2$, and the image is an Enriques sextic. These last surfaces form an irreducible and unirational family of dimension 6 of the moduli space.*

*Moreover,*

- *if $K^2 = 8$, $S$ has an isotrivial genus 3 rational fibration whose general fibre is hyperelliptic with 6 double fibres;*
- *if $K^2 = 7$, $S$ has a genus 3 rational fibration whose general fibre is hyperelliptic with 5 double fibres and a fibre with reducible support, consisting of two components;*
- *if $K^2 = 6$ and $\deg \varphi_{2K} = 2$, $S$ has a genus 3 rational fibration whose general fibre is hyperelliptic with 4 or 5 double fibres;*
- *if $K^2 = 6$ and $\deg \varphi_{2K} = 4$, $S$ is a Burniat surface.*

*Remark 3.* Surfaces with $p_g = q = 0$ and $K^2 = 3, 4$ were constructed by Keum ([Keu88]) and Naie ([Nai94]). For $K^2 = 3$ the degree of the bicanonical map can be equal to 2 and to 4, and it is an open question if it can be birational.

Concerning the classification of surfaces with $p_g = q = 0$ there has been recent progress.

We would like to mention that in [BCG05b] a complete classification of surfaces with $p_g = q = 0$ isogenous to a product (this forces $K^2 = 8$) is given.

In the case $p_g = q = 0$ and $K^2 = 9$ all the surfaces in question are, by Yau's theorem ([Yau77] and [Yau78]) quotient of the complex unit ball in $\mathbb{C}^2$ by a discrete group $\Gamma$ acting freely. The first effective example of such surfaces, called *fake projective planes* since they have the same Betti numbers as the projective plane, was given by Mumford ([Mum79]) using 2-adic uniformization. Other examples were given later in [IK98], while recently Keum ([Keu05]) gave an explicit geometric construction as a cyclic cover of a particular Dolgachev surface.

In a recent preprint ( [PY05]) G. Prasad and S.K. Yeung, using the arithmeticity of $\Gamma$, show that the fake projective planes belong to twelve rather explicit lists, each corresponding to an imaginary quadratic field $\mathbb{Q}(\sqrt{-a})$ and a prime $p$ which ramifies in it.

The interesting geometric corollaries are that:

i) all these groups $\Gamma$ are indeed contained in $SU(2,1)$, hence the canonical divisor $K$ is divisible by 3,

ii) all these surfaces have a nontrivial first homology group $H_1(S, \mathbb{Z})$.

## 3 Surfaces with $p_g = 4$

In Enriques' book on algebraic surfaces [Enr49] much emphasis was put on the effective construction of surfaces whose canonical map is birational, particularly for surfaces with $p_g = 4$, where the canonical image is a surface in $\mathbb{P}^3$.

Later on, in particular in the last thirty years, many authors studied surfaces with $p_g = 4$, with particular interest in the construction of surfaces with $p_g = 4$, birational canonical map and $K^2$ as high as possible.

If the canonical map of a minimal surface of general type $S$ with $p_g = 4$ is birational, then the standard inequalities give $5 \leq K^2 \leq 45$.

Nowadays we know examples, by the contribution of several authors, for every value of $K_S^2$ in the range $5 \leq K_S^2 \leq 28$ (cf. e.g. [Cil81], [Cat99]). An example with $K_S^2 = 31$ has been recently obtained in [Lie03], although the example constructed has a big fixed part of the canonical system so that its canonical image has "only" degree 12. Moreover, the first two authors together with F. Grunewald have constructed a canonical surface in the projective 3-space with $K^2 = 45$. This surface is obtained as a Galois covering of the plane with group $(\mathbb{Z}/5\mathbb{Z})^2$, branched over a configuration of lines introduced by Hirzebruch (cf. [BCG05c]).

In this case we have a rigid surface such that its canonical system has a fixed part.

Obviously also in this case classification is the biggest challenge: for which values of $K_S^2$ is it possible to classify all possible minimal surfaces of general type with $p_g = 4$? And more ambitiously: for which values of $K_S^2, q$ it is possible to completely describe the moduli space $\mathcal{M}_{K^2, 4, q}$?

### 3.1 $K^2 = 4, 5$

The cases $K^2 = 4, 5$ were already treated by Enriques ([Enr49], section 2, chapter VIII, pp.268–271), and the corresponding moduli spaces were completely understood already in the 70's.

We briefly recall these results.

By Debarre's inequality, all these surfaces are regular (in fact, this is true for $K^2 \leq 7$). The canonical map of surfaces with $K^2 = 4$ and $p_g = 4$ is a morphism of degree 2 onto an irreducible quadric in $\mathbb{P}^3$, i.e., either a smooth quadric or a quadric cone. The general surface is a double cover of a smooth quadric branched over a smooth complete intersection with a sextic surface.

A detailed analysis of the corresponding moduli space can be found in [Hor76], where the following is proven.

**Theorem 17.** $\mathcal{M}_{4,4,0}$ *is irreducible, unirational of dimension* 42, *its singular locus is irreducible of codimension* 1 *and corresponds exactly to the surfaces whose canonical image is a quadric cone.*

We have two classes of minimal surfaces of general type with $K^2 = 4$: let us say surfaces of type I (double covers of a smooth quadric) and II (double covers of a quadric cone). Correspondingly we have a stratification of $\mathcal{M}_{4,4,0}$ as a union of two locally closed strata, both irreducible, which we denote simply by I and II, of respective dimension 42 and 41.

To draw a picture of this moduli space we need the following notation:

**Definition 9.** *Let A and B be two (locally closed) irreducible strata of a moduli space $\mathcal{M}_{K^2,p_g,q}$.*

*If we write "$A \to B$", it means that there is a flat family with base a small disc $\Delta_\varepsilon \subset \mathbb{C}$, whose central fibre is of type B and whose general fibre is of type A. In other words it means that the closure of the stratum A intersects the stratum B.*

With this notation a picture of $\mathcal{M}_{4,4,0}$ is the following:

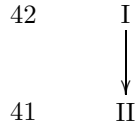

Note that at the left of each stratum stands the dimension of the corresponding irreducible stratum.

The case $K^2 = 5$ is slightly more complicated, and completely described in [Hor75]: the canonical map is either a birational morphism to a quintic in $\mathbb{P}^3$ (type I), or a rational map of degree 2 onto an irreducible quadric, which can be as in the previous case either smooth (type $II_a$) or a quadric cone (type $II_b$).

**Theorem 18.** $\mathcal{M}_{5,4,0}$ *has two irreducible components, both unirational of dimension* 40, *intersecting in a* 39 *dimensional subvariety.*

Here is the picture for $\mathcal{M}_{5,4,0}$

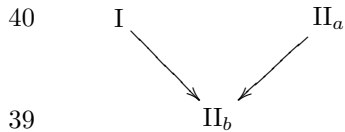

The moduli space has two irreducible unirational components of dimension 40 whose general point corresponds to surfaces with canonical image respectively a quintic or a smooth quadric. The surfaces whose canonical image is a quadric cone form a 39−dimensional subvariety of this moduli space, the intersection of the two irreducible components.

## 3.2 $K^2 = 6$

This case is much more complicated and surfaces with $K^2 = 6$ and $p_g = 4$ were completely classified by Horikawa in [Hor78], obtaining a stratification of $\mathcal{M}_{6,4,0}$ in 11 strata. We will not enter here the details of this classification.

A complete understanding of $\mathcal{M}_{6,4,0}$ is still missing, since it is not clear how exactly these 11 strata "glue" toghether.

**Theorem 19 ([Hor78],[BCP04]).** $\mathcal{M}_{6,4,0}$ *has* 4 *irreducible components, all unirational, one of dimension* 39, *the others of dimension* 38.

*Moreover, the number of connected components of* $\mathcal{M}_{6,4,0}$ *is at most two.*

The main question on this moduli space remains the following:

*Question 1.* Is $\mathcal{M}_{6,4,0}$ connected?

Here is a partial picture:

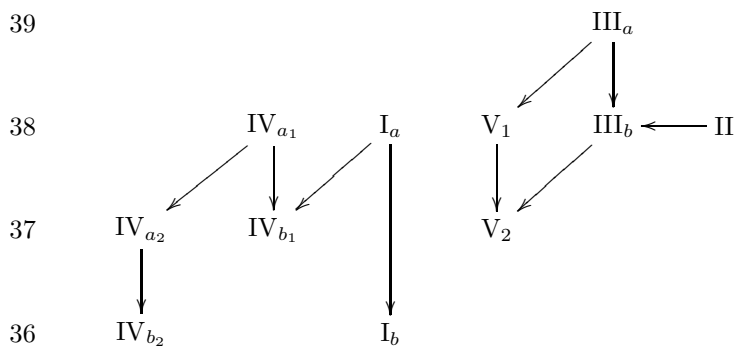

This picture is partial because up to now it is not known whether all possible arrows are drawn. More precisely, $\mathcal{M}_{6,4,0}$ is connected if and only if one of the two following degenerations is possible: $I_a \to V_1$ or $I_a \to V_2$.

This picture was done by Horikawa in [Hor78] with the exception of the horizontal line $III_b \leftarrow II$, recently obtained in [BCP04].

We are going to explain how this arrow was obtained.

We need to construct a flat family of surfaces whose central fibre is of type $III_b$ and whose general fibre is of type II. The difficulty lies in the obvious fact that, by reasons of dimension, the general surface of type $III_b$ cannot deform to a surface of type II.

First however we need to explain what surfaces of type II and of type III$_b$ are.

Surfaces of type II are defined as follows

**Definition 10.** *A minimal surface of general type with $p_g = 4$ and $K^2 = 6$ is of type II if the canonical map has degree* 3.

It is immediate, that the canonical image is a quadric cone.

In fact, one can see that their canonical models $X$ are hypersurfaces of degree 9 in the weighted projective space $\mathbb{P}(1,1,2,3)$, so that the canonical divisor of $X$ is divisible by 2 as a Weil divisor.

**Definition 11.** *A minimal surface of general type with $p_g = 4$ and $K^2 = 6$ is of type* III$_b$ *if it has no genus* 2 *pencil and the canonical system has a fixed component.*

Horikawa gave a very concrete description of this class (theorem 5.2 in [Hor78]) showing, among other things, that the canonical map is a double cover of a quadric cone.

In both cases there is a pencil $L$ on $S$, the strict transform of the ruling of the quadric cone, such that $K_S - 2L$ is effective.

Computing intersection numbers one sees (cf. also [MLP00] for case II) that in both cases we have a decomposition $K_S = 2L + Z$ where $Z$ is a fundamental cycle ($Z^2 = -2$, $KZ = 0$) and $L$ is a genus 3 pencil with one simple base point. The main difference between the two cases is that in case III$_b$ all curves in the pencil are hyperelliptic, while in case II the general curve is nonhyperelliptic.

The idea is then to construct a family preserving this decomposition and therefore the genus 3 pencil. Since $Z$ is a fundamental cycle, one can consider the canonical model $X$ of $S$, the surface with rational double points obtained contracting all fundamental cycles. The canonical class of $X$ is then 2−divisible (as a Weil divisor), since $K_X = 2\tilde{L}$, $\tilde{L}$ being the image of $L$ on $X$.

We consider the semicanonical ring $R = R(X, \tilde{L})$: the name *semicanonical* being selfexplanatory, since the subring generated by the homogeneous elements of even degree is exactly the canonical ring of $X$ (and of $S$).

In order to compute the ring $R$ we use the hyperplane section principle [Rei90a], first computing the quotient ring $R(X, \tilde{L})/(x_0)$, where $x_0$ is a general homogeneous element of degree 1, i.e., corresponding to a general element $C$ of the pencil $\tilde{L}$.

It is not difficult to prove that this ring equals the ring

$$R(C, \frac{3}{2}P) := \bigoplus_m H^0(C, \mathcal{O}_C([\frac{3m}{2}]P))$$

(we refer for the definition of the ring structure to [BCP04]), which is a suitable subring of the ring $R(C, P)$, where $P$ is the base point of $L$ (and a Weierstrass point for the hyperelliptic curve $C$).

Then the strategy is the following:

i) take any hyperelliptic genus 3 curve $C$, a Weierstraß point $P$ on it, and compute the ring $R(C, \frac{3}{2}P)$
ii) "deform" it, adding an element of degree 1 in a flat way, to get the semi-canonical ring of a surface of type $III_b$
iii) construct a flat family of rings (say with parameter $t$), whose central fibre ($t = 0$) is the ring constructed in ii), and whose general fibre ($t \neq 0$) is a hypersurface ring of type II

Step i) is easy, steps ii) and iii) depend on the result of step i), concerning $R(C, \frac{3}{2}P)$.

It turns out that $R(C, \frac{3}{2}P)$ is a Gorenstein ring of codimension 4, which can be expressed in a nice way:

**Proposition 1.** *Let $C$ be a hyperelliptic curve of genus 3, $p \in C$ a Weierstraß point. Then $R(C, \frac{3}{2}p) \cong \mathbb{C}[x, y, z, w, v, u]/I$, where $\deg(x, y, z, w, v, u) = (1, 2, 3, 4, 5, 6)$ and the ideal $I$ is generated by the $4 \times 4$ Pfaffians of the skew-symmetric 'extra-symmetric' matrix*

$$M = \begin{pmatrix} 0 & 0 & z & v & y & x \\ & 0 & w & u & z & y \\ & & 0 & \tilde{P}_9 & u & v \\ & & & 0 & w^2 & zw \\ & & & & 0 & 0 \\ -sym & & & & & 0 \end{pmatrix},$$

*where $\tilde{P}_9$ is homogeneous of degree 9 in the variables $x, y, z, w$.*

The graded matrix $M$ has the nice property to be "extrasymmetric".
*Extrasymmetric matrices.* Extrasymmetric matrices were introduced by Miles Reid and Duncan Dicks ([Rei90a],[Rei89]). Let $A$ be a polynomial ring and let $M$ be a skew 'extrasymmetric' matrix of the form

$$\begin{pmatrix} 0 & a & b & c & d & e \\ & 0 & f & g & h & d \\ & & 0 & i & g & c \\ & & & 0 & pf & pb \\ & & & & 0 & pa \\ -sym & & & & & 0 \end{pmatrix}$$

where $a, b, c, d, e, f, g, h, i, p \in A$. Then 6 of the 15 ($4 \times 4$) Pfaffians belong to the ideal generated by the other 9; moreover, if the entries are general enough, the ideal generated by these pfaffians has exactly 16 independent syzygies, *which can all be written explicitly as functions of the entries of the matrix.*

This implies that, if we have a ring presented in this form, and the ring has no further syzygies, deforming the entries of the matrix (preserving the symmetries), we obtain automatically a a flat deformation of the ring. We mention that recent studies lead to generalizations of this "format": here we need only this special "original" case.

By proposition 1 $R(C, \frac{3}{2}P)$ is presented by an extrasymmetric matrix. Therefore we are reduced to add "$x_0$" (step ii)) and "$t$" (step iii)) such that the obtained matrix is still extrasymmetric and has homogeneous entries and Pfaffians. Moreover, we have to take care that for $t = 0$ we obtain a fibration by hyperelliptic curves, whereas for $t \neq 0$ the general curve of the fibration has to be non hyperelliptic.

A crucial ingredient is

*The zero of degree zero.* Notice that the second entry of the first row of $M$ is equal to 0, and it corresponds to a homogeneous element of degree 0. This allows to substitute in this entry a parameter $t$. For $t \neq 0$ the upshot is that from the 9 Pfaffians we can eliminate the variables $w, v, u$ and we are left with the variables $x_0, x_1, y, z$ and with a single equation of degree 9: that is, we have a semicanonical ring of type II, and we have succeeded!

We obtain thus the following result:

**Theorem 20.** *Consider the ring $\mathbb{C}[x_0, x_1, y, z, w, v, u]$ with variables of respective degrees $(1, 1, 2, 3, 4, 5, 6)$.*

*Consider a family of skew extrasymmetric matrices, with parameter $t$*

$$M_t = \begin{pmatrix} 0 & t\,z & v & y & x_1 \\ & 0 & w & u & P_3 & y \\ & & 0 & P_9 & u & v \\ & & & 0 & wP_4 & zP_4 \\ & & & & 0 & tP_4 \\ -sym & & & & & 0 \end{pmatrix}.$$

*where the $P_i$'s are homogeneous of degree $i$ in the first 5 variables of the ring, and let $J_t$ be the ideal generated by the $4 \times 4$ pfaffians of $M_t$.*

*Then, for general choice of the polynomials $P_i$, $\mathbb{C}[x_0, x_1, y, z, w, v, u]/J_t$ is, for $t = 0$, the semicanonical ring of a surface of type $\mathrm{III}_b$, and for $t \neq 0$ the semicanonical ring of a surface of type II.*

The surfaces of type $\mathrm{III}_b$ whose semicanonical ring can be presented as in the above theorem form a codimension 2 subscheme of the corresponding 38−dimensional stratum of $\mathcal{M}_{6,4,0}$, lying in the intersection with the 38−dimensional stratum II.

It is still unclear whether this is the whole intersection of the closures of these two strata. A priori one can only say that this intersection has at most dimension 37, and what we found has only dimension 36. Therefore there remains the following:

*Question 2.* Exactly which surfaces of type $\mathrm{III}_b$ lie in the closure of the stratum II?

The above question is of course related to the singularity type of the local moduli space.

### 3.3 $K^2 = 7$

This is the last case for which there is a complete classification.

In [Bau01] the first author gives a very precise description of surfaces with $p_g = 4$, $K^2 = 7$ according to the behaviour of the canonical map, allowing to show that the moduli space $\mathcal{M}_{7,4,0}$ has three irreducible components $\mathcal{M}_{36}$, $\mathcal{M}'_{36}$ and $\mathcal{M}_{38}$ of respective dimensions 36, 36 and 38. Moreover, it is shown that the two irreducible components of dimension 36 intersect, whereas it is up to now not yet clear whether the component of dimension 38 is indeed a connected component or it intersects $\mathcal{M}_{36}$.

We encounter here a very similar situation as for $K^2 = 6$. There are two families, one in $\mathcal{M}_{36}$, the other in $\mathcal{M}_{38}$, where the first consists of surfaces admitting a non hyperelliptic genus 3 pencil, whereas the surfaces in the other family admit a hyperelliptic genus 3 pencil. In fact, the first and last author have been able to calculate the relative canonical algebra for the hyperelliptic case, which is Gorenstein of codimension 6. For this very high codimension there are yet no flexible formats known to organize the equations. Hopefully it will be possible to understand the deformations of this family.

## 4 Surfaces isogeneous to a product, Beauville surfaces and the absolute Galois group

Surfaces isogenous to a (higher) product were introduced and extensively studied by the second author in [Cat00], where it is proven that any surface $S$ isogenous to a higher product has a unique minimal realization as a quotient $S = (C_1 \times C_2)/G$.

Here $C_1$ and $C_2$ are smooth algebraic curves of genus at least 2 and $G$ is a finite group acting freely, and with the property that no element acts trivially on one of the factors $C_i$.

Moreover, it was shown that the topology of a surface isogenous to a product determines its deformation class up to complex conjugation. The following result contains a correction to Theorem 4.14 of [Cat00] (cf. theorem 3.3 of [Cat03]).

**Theorem 21.** *Let $S = (C_1 \times C_2)/G$ be a surface isogenous to a product. Then any surface $S'$ with the same topological Euler number and the same fundamental group as $S$ is diffeomorphic to $S$. If moreover $S'$ is orientedly diffeomorphic to $S$, then $S'$ is deformation equivalent to $S$ or to $\bar{S}$. In other words, the corresponding moduli space $\mathcal{M}_S^{top} = \mathcal{M}_S^{diff}$ is either irreducible and connected or it contains two connected components which are exchanged by complex conjugation.*

This class of surfaces and their higher dimensional analogues provide a wide specimen of examples where one can test or disprove several conjectures and questions (cf. e.g. [Cat03], [BC04], [BCG05a], compare also the next section).

Moreover, the absolute Galois group $Aut(\bar{\mathbb{Q}}/\mathbb{Q})$ acts on the moduli spaces of this class of surfaces: we shall outline a direct connection with Grothendieck's dream of "dessins d'enfants".

In the following we shall concentrate on a class of surfaces isogenous to a product, namely the rigid ones.

We recall that an algebraic variety $X$ is *rigid* if and only if it does not have any non trivial deformations (e.g., the projective space is rigid). There is another (stronger) notion of rigidity, which is the following

**Definition 12.** *An algebraic variety $X$ is called* strongly rigid *if any other variety homotopically equivalent to $X$ is either biholomorphic or antibiholomorphic to $X$.*

*Remark 4.* 1) It is nowadays wellknown that smooth compact quotients of symmetric spaces are rigid (cf. [CV60]).

2) Mostow (cf. [Mos73] proved that indeed locally symmetric spaces of complex dimension $\geq 2$ are strongly rigid, in the sense that any homotopy equivalence is induced by a unique isometry.

These varieties are of general type and the moduli space of varieties of general type is defined over $\mathbb{Z}$, and naturally the absolute Galois group $Gal(\bar{\mathbb{Q}}/\mathbb{Q})$ acts on the set of their connected components. So, in our special case, $Gal(\bar{\mathbb{Q}}/\mathbb{Q})$ acts on the isolated points which parametrize rigid varieties.

In particular, rigid varieties are defined over a number field and work of Shimura gives a possible way of computing explicitly their fields of definition. By this reason these varieties were named *Shimura varieties* (cf. Deligne's Bourbaki seminar [Del71]).

A quite general question is

*Question 3.* What are the fields of definition of rigid varieties? What is the $Gal(\bar{\mathbb{Q}}/\mathbb{Q})$-orbit of the point in the moduli space corresponding to a rigid variety?

Much simpler examples of rigid varieties were found by the second author (cf. [Cat00]).

**Beauville surfaces**

Inspired by a construction of A. Beauville of a surface with $K^2 = 8, p_g = q = 0$ (cf. [Bea78]) as a quotient of the product of two Fermat curves of degree 5 by the action of the group $\mathbb{Z}/5\mathbb{Z}$, in [Cat00] the following definition was given

**Definition 13.** *A* Beauville *surface is a compact complex surface $S$ which*
  *1) is* rigid*, i.e., it has no nontrivial deformation,*

*2) is* isogenous to a higher product, *i.e., it is a quotient* $S = (C_1 \times C_2)/G$ *of a product of curves of resp. genera* $\geq 2$ *by the free action of a finite group* $G$.

Notice that, given a surface isogenous to a product, we obtain always three more, exchanging $C_1$ with its conjugate curve $\bar{C}_1$, or $C_2$ with $\bar{C}_2$: but only if we conjugate both $C_1$, $C_2$ we obtain an orientedly diffeomorphic surface. These four surfaces could however be all biholomorphic to each other.

If $S$ is a Beauville surface and $X$ is orientedly diffeomorphic to $S$, then theorem 21 implies: $X \cong S$ or $X \cong \bar{S}$.

In other words, the corresponding subset of the moduli space $\mathfrak{M}_S$ consists of one or two points (if we insist on keeping the orientation fixed, else we may get up to four points).

**Definition 14.** *C is a* triangle curve *if there is a finite group G acting effectively on C and satisfying the properties*
 *i)* $C/G \cong \mathbb{P}^1_{\mathbb{C}}$, *and*
 *ii)* $f : C \to \mathbb{P}^1_{\mathbb{C}} \cong C/G$ *has* $\{0, 1, \infty\}$ *as branch set.*

*Remark 5.* The rigidity of a Beauville surface is equivalent to the condition that $(C_i, G^0)$ is a triangle curve, for $i = 1, 2$ ($G^0 \subset G$ is the subgroup of index $\leq 2$ which does not exchange the two factors).

Recall now the classical

**Theorem 22.** (Riemann's Existence Theorem)
 *There is a natural bijection between:*
 *1) Equivalence classes of holomorphic mappings* $f : C \to \mathbb{P}^1_{\mathbb{C}}$, *of degree n and with Branch set* $B_f \subset B$, *(where C is a compact Riemann surface, and* $f : C \to \mathbb{P}^1_{\mathbb{C}}$, $f' : C' \to \mathbb{P}^1_{\mathbb{C}}$ *are said to be equivalent if there is a biholomorphism* $g : C' \to C$ *such that* $f' = f \circ g$).
 *2) Conjugacy classes of* monodromy *homomorphisms* $\mu : \pi_1(\mathbb{P}^1_{\mathbb{C}} - B) \to \mathfrak{S}_n$ *(here,* $\mathfrak{S}_n$ *is the symmetric group in n letters, and* $\mu \cong \mu'$ *iff there is an element* $\tau \in \mathfrak{S}_n$ *with* $\mu(\gamma) = \tau \mu'(\gamma) \tau^{-1}$, $(\forall \gamma)$.
 *Moreover:*
 *3) C is connected if and only if the subgroup* $Im(\mu)$ *acts transitively on* $\{1, 2, \ldots n\}$.
 *4) f is a polynomial if and only if* $\infty \in B$, *the monodromy at* $\infty$ *is a cyclical permutation, and* $g(C) = 0$.

*Remark 6.* 1) Assume that $\infty \in B$, so $\{\infty, b_1, \ldots b_d\} = B$: then $\pi_1(\mathbb{P}^1_{\mathbb{C}} - B)$ is a free group generated by $\gamma_1, \ldots \gamma_d$ and $\mu$ is completely determined by the local monodromies $\tau_i := \mu(\gamma_i)$.

Grothendieck's enthusiasm was raised by the following result, where Belyi ([Bel79])made a very clever and very simple use of some explicit polynomials, now called the Belyi polynomials, of the form $\frac{(m+r)^{m+r}}{m^m r^r} z^m (z-1)^r$ in order to reduce the number of critical values of an algebraic function defined over $\bar{\mathbb{Q}}$.

**Theorem 23.** (Belyi) *An algebraic curve $C$ can be defined over $\overline{\mathbb{Q}}$ if and only if there exists a holomorphic map $f : C \to \mathbb{P}^1_{\mathbb{C}}$ with branch set only $\{0, 1, \infty\}$.*

The word "dessin d' enfant" = child's drawing is due to the fact that the monodromy of $f$ is determined by the 'dessin d' enfant' $f^{-1}([0, 1])$, a bipartite graph (the vertices have label 0 or 1 according to their image) such that at each vertex one has a cyclical order of the edges incident in the vertex (this property holds because the graph is contained in a complex curve $C$ thus we choose the corresponding cyclical counterclockwise order).

It is clear that the *triangle curves* correspond to a certain class of 'dessins d' enfants', those which admit a group action with quotient the interval $[0, 1]$.

Let us parenthetically observe that Gabino Gonzalez was recently able to extend Belyi's theorem to the case of complex surfaces (in terms of Lefschetz maps with three critical values) (cf. [Gon04]).

Grothendieck ([Gro97] )proposed to look at the 'dessins d' enfants' in order to get representations of the absolute Galois group $Gal(\overline{\mathbb{Q}}, \mathbb{Q})$.

We just explained that a Beauville surface is defined over $\overline{\mathbb{Q}}$, and that the Galois group $Gal(\overline{\mathbb{Q}}, \mathbb{Q})$ operates on the discrete subset of the moduli space $\mathfrak{M}_S$ corresponding to Beauville surfaces.

This action is rather strictly related to the action on the 'dessins d' enfants', but in this case, by theorem 21, the Galois group $Gal(\overline{\mathbb{Q}}, \mathbb{Q})$ may transform a Beauville surface into another one with a non isomorphic fundamental group.

Phenomena of this kind were already observed by J.P. Serre (cf. [Ser64]): here the idea is not to consider this as a pathology, but as a source of information, and to actually try to understand the representation of the Galois group $Gal(\overline{\mathbb{Q}}, \mathbb{Q})$ on the class of groups which are fundamental groups of Beauville surfaces (and of their higher dimensional analogues).

It looks therefore interesting to investigate these surfaces and to address the following problems:

*Question 4.* Existence and classification of Beauville surfaces, i.e.,
   a) which finite groups $G$ can occur?
   b) classify all possible Beauville surfaces for a given finite group $G$.

*Question 5.* Is the Beauville surface $S$ biholomorphic to its complex conjugate surface $\bar{S}$?

Is $S$ real (i.e., does there exist a biholomorphic map $\sigma : S \to \bar{S}$ with $\sigma^2 = id$)?

Another motivation to find these surfaces was also given by the following
**FRIEDMAN-MORGAN'S SPECULATION ( [FM88] 1987):**
DEF $\iff$ DIFF (Differentiable equivalence and deformation equivalence coincide for surfaces).

A series of counterexamples were given by several authors and the simplest examples were given using non rigid surfaces isogenous to a product. Great part of the next section will be devoted to these equivalence relations.

In order to reduce the description of Beauville surfaces to some group theoretic statement, we need to recall that surfaces isogenous to a higher product belong to two types:

- $S$ is of *unmixed type* if the action of $G$ does not mix the two factors, i.e., it is the product action of respective actions of $G$ on $C_1$, resp. $C_2$.
- $S$ is of *mixed type*, i.e., $C_1$ is isomorphic to $C_2$, and the subgroup $G^0$ of transformations in $G$ which do not mix the factors has index precisely 2 in $G$.

The datum of a Beauville surface can be completely described group theoretically, since it is equivalent to the datum of two triangle curves with isomorphic groups.

**Definition 15.** *Let $G$ be a finite group.*
*1) A quadruple $v = (a_1, c_1; a_2, c_2)$ of elements of $G$ is an* unmixed Beauville *structure for $G$ if and only if*
*(i) the pairs $a_1, c_1$, and $a_2, c_2$ both generate $G$,*
*(ii) $\Sigma(a_1, c_1) \cap \Sigma(a_2, c_2) = \{1_G\}$, where*

$$\Sigma(a, c) := \bigcup_{g \in G} \bigcup_{i=0}^{\infty} \{g a^i g^{-1}, g c^i g^{-1}, g(ac)^i g^{-1}\}.$$

*We write $\mathbb{U}(G)$ for the set of unmixed Beauville structures on $G$.*
*2) A mixed Beauville quadruple for $G$ is a quadruple $M = (G^0; a, c; g)$ consisting of a subgroup $G^0$ of index 2 in $G$, of elements $a, c \in G^0$ and of an element $g \in G$ such that*
*i) $G^0$ is generated by $a, c$,*
*ii) $g \notin G^0$,*
*iii) for every $\gamma \in G^0$ we have $g\gamma g\gamma \notin \Sigma(a, c)$.*
*iv) $\Sigma(a, c) \cap \Sigma(gag^{-1}, gcg^{-1}) = \{1_G\}$.*
*We write $\mathbb{M}(G)$ for the set of mixed Beauville quadruples on the group $G$.*

*Remark 7.* We consider here finite groups $G$ having a pair $(a, c)$ of generators. Setting $(r, s, t) := (\text{ord}(a), \text{ord}(c), \text{ord}(ac))$, such a group is a quotient of the triangle group

$$T(r, s, t) := \langle x, y \mid x^r = y^s = (xy)^t = 1 \rangle. \tag{1}$$

It is now easy to explain how to get a surface from the above data, if we remember Riemann's Existence Theorem which we recalled just above.

We take as base point $\infty \in \mathbb{P}^1_{\mathbb{C}}$ and consider $B := \{-1, 0, 1\}$. We choose the following generators $\alpha, \beta$ of $\pi_1(\mathbb{P}^1_{\mathbb{C}} - B, \infty)$ ($\gamma := (\alpha \cdot \beta)^{-1}$):

Let now $G$ be a finite group and $v = (a_1, c_1; a_2, c_2) \in \mathbb{U}(G)$. We get surjective homomorphisms

$$\pi_1(\mathbb{P}^1_{\mathbb{C}} - B, \infty) \to G, \qquad \alpha \mapsto a_i, \ \gamma \mapsto c_i \tag{2}$$

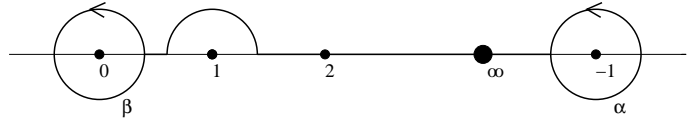

and Galois coverings $\lambda_i : C(a_i, c_i) \to \mathbb{P}^1_{\mathbb{C}}$ ramified only in $\{-1, 0, 1\}$ with ramification indices equal to the orders of $a_i$, $b_i$, $c_i$ and with group $G$ (by *Riemann's existence theorem*).

*Remark 8.* 1) Condition (1), ii) ensures that the action of $G$ on $C(a_1, c_1) \times C(a_2, c_2)$ is free.

2) Let be $\iota(a_1, c_1; a_2, c_2) = (a_1^{-1}, c_1^{-1}; a_2^{-1}, c_2^{-1})$. Then $S(\iota(v)) = \overline{S(v)}$ (note in fact that $\bar{\alpha} = \alpha^{-1}$, $\bar{\gamma} = \gamma^{-1}$).

3) One can verify that the required conditions automatically imply: $g(C(a_1, c_1)) \geq 2$ and $g(C(a_2, c_2)) \geq 2$.

One has:

**Proposition 2.** *Let $G$ be a finite group and*

$v = (a_1, c_1; a_2, c_2) \in \mathbb{U}(G)$.

*Assume that* $\{\text{ord}(a_1), \text{ord}(c_1), \text{ord}(a_1c_1)\} \neq \{\text{ord}(a_2), \text{ord}(c_2), \text{ord}(a_2c_2)\}$ *and that* $\text{ord}(a_i) < \text{ord}(a_ic_i) < \text{ord}(c_i)$. *Then $S(v) \cong \overline{S(v)}$ if and only if there are inner automorphisms $\phi_1, \phi_2$ of $G$ and an automorphism $\psi \in Aut(G)$ such that, setting $\psi_j := \psi \circ \phi_j$, we have $\psi_1(a_1) = a_1^{-1}$, $\psi_1(c_1) = c_1^{-1}$, and $\psi_2(a_2) = a_2^{-1}$, $\psi_2(c_2) = c_2^{-1}$.*

*In particular, under the above assumption, $S(v)$ is isomorphic to $\overline{S(v)}$ if and only if $S(v)$ has a real structure.*

*Remark 9.* Dropping the assumption on the orders of $a_i$, $c_i$, we can define a finite permutation group $A_{\mathbb{U}}(G)$ such that for $v, v' \in \mathbb{U}(G)$ we have: $S(v) \cong S(v')$ if and only if $v$ is in the $A_{\mathbb{U}}(G)$-orbit of $v'$.

*Remark 10.* If $G$ is abelian, $v \in \mathbb{U}(G)$. Then $S(v)$ always has a real structure.

We have the following results (cf. [BCG05a] for some of these, others have not yet been published):

**Theorem 24.** *1) An abelian group $G$ admits an unmixed Beauville structure iff $G \cong (\mathbb{Z}/n)^2$, $(n, 6) = 1$.*

*2) The following groups admit unmixed Beauville structures:*

*a) the alternating group $\mathfrak{A}_n$ for large $n$,*

*b) the symmetric group $\mathfrak{S}_n$ for $n \in \mathbb{N}$ with $n \geq 7$*

*c) $\mathbf{SL}(2, \mathcal{F}_p)$, $\mathbf{PSL}(2, \mathcal{F}_p)$ for $p \neq 2, 3, 5$.*

With the help of the computer algebra program MAGMA all finite simple nonabelian groups of order $\leq 50000$ were checked and unmixed Beaville structures were found on all of them, with the exception of $\mathfrak{A}_5$ (where it can't be found!). This led to the following

*Conjecture 1.* ([BCG05a]) All finite simple nonabelian groups except $\mathfrak{A}_5$ admit an unmixed Beauville structure.

This conjecture was also checked for some bigger simple groups like the Mathieu groups **M**12, **M**22 and also matrix groups of size bigger then 2.

Call now $(r, s, t) \in \mathbb{N}^3$ *hyperbolic* if

$$\frac{1}{r} + \frac{1}{s} + \frac{1}{t} < 1.$$

In this case the triangle group $T(r, s, t)$ is hyperbolic. These studies led also to the following suggestive:

*Conjecture 2.* ([BCG05a]) Let $(r, s, t)$, $(r', s', t')$ be two hyperbolic types. Then almost all alternating groups $\mathfrak{A}_n$ have an unmixed Beauville structure $v = (a_1, c_1; a_2, c_2)$ where $(a_1, c_1)$ has type $(r, s, t)$ and $(a_2, c_2)$ has type $(r', s', t')$.

The above conjectures are variations of a conjecture of Higman (proved by B. Everitt (2000), [Eve00]) asserting that every hyperbolic triangle group surjects onto almost all alternating groups.

Concrete explicit examples of rigid surfaces not biholomorphic to their complex conjugate were also given:

**Theorem 25.** *The following groups admit unmixed Beauville structures $v$ such that $S(v)$ is not biholomorpic to $\overline{S(v)}$:*

*1. the symmetric group $\mathfrak{S}_n$ for $n \geq 7$*

*2. the alternating group $\mathfrak{A}_n$ for $n \geq 16$ and $n \equiv 0 \mod 4$, $n \equiv 1 \mod 3$, $n \not\equiv 3, 4 \mod 7$.*

And also new examples of real points of moduli spaces which do not correspond to real surfaces:

**Theorem 26.** *Let $p > 5$ be a prime with $p \equiv 1 \mod 4$, $p \not\equiv 2, 4 \mod 5$, $p \not\equiv 5 \mod 13$ and $p \not\equiv 4 \mod 11$. Set $n := 3p + 1$. Then there is an unmixed Beauville surface $S$ with group $\mathfrak{A}_n$ which is biholomorphic to the complex conjugate surface $\bar{S}$, but is not real.*

For *mixed* Beauville surfaces the situation is more complicated, as already the following suggests.

**Theorem 27.** *1) If a group $G$ admits a mixed Beauville structure, then the subgroup $G^0$ is non abelian.*

*2) No group of order $\leq 512$ admits a mixed Beauville structure.*

A general construction of finite groups admitting a mixed Beauville structure was given in [BCG05a].

Let $H$ be a non-trivial group, and let $\Theta : H \times H \to H \times H$ be the automorphism defined by $\Theta(g, h) := (h, g)$ $(g, h \in H)$. We consider the semidirect product

$$H_{[4]} := (H \times H) \rtimes \mathbb{Z}/4\mathbb{Z} \qquad (3)$$

where the generator 1 of $\mathbb{Z}/4\mathbb{Z}$ acts through $\Theta$ on $H \times H$. Since $\Theta^2$ is the identity we find

$$H_{[2]} := H \times H \times (2\mathbb{Z}/4\mathbb{Z}) \cong H \times H \times \mathbb{Z}/2\mathbb{Z} \qquad (4)$$

as a subgroup of index 2 in $H_{[4]}$.

We have now

**Lemma 1.** *Let $H$ be a non-trivial group and let $a_1, c_1, a_2, c_2$ be elements of $H$. Assume that*

1. *the orders of $a_1, c_1$ are even,*
2. $a_1^2, a_1 c_1, c_1^2$ *generate $H$,*
3. $a_2, c_2$ *also generate $H$,*
4. $(\mathrm{ord}(a_1) \cdot \mathrm{ord}(c_1) \cdot \mathrm{ord}(a_1 c_1), \mathrm{ord}(a_2) \cdot \mathrm{ord}(c_2) \cdot \mathrm{ord}(a_2 c_2)) = 1$.

*Set $G := H_{[4]}$, $G^0 := H_{[2]}$ as above and $a := (a_1, a_2, 2)$, $c := (c_1, c_2, 2)$. Then $(G^0; a, c)$ is a mixed Beauville structure on $G$.*

*Proof.* It is easy to see that $a, c$ generate $G^0 := H_{[2]}$.

The crucial observation is that

$$(1_H, 1_H, 2) \notin \Sigma(a, c). \qquad (5)$$

In fact, if this were not correct, it would have to be conjugate of a power of $a$, $c$ or $b$. Since the orders of $a_1$, $b_1$, $c_1$ are even, we obtain a contradiction.

Suppose that $h = (x, y, z) \in \Sigma(a, c)$ satisfies $\mathrm{ord}(x) = \mathrm{ord}(y)$: then our condition 4 implies that $x = y = 1_H$ and (5) shows $h = 1_{H_{[4]}}$.

Let now $g \in H_{[4]}$, $g \notin H_{[2]}$ and $\gamma \in G^0 = H_{[2]}$ be given. Then $g\gamma = (x, y, \pm 1)$ for appropriate $x, y \in H$. We find

$$(g\gamma)^2 = (xy, yx, 2)$$

and the orders of the first two components of $(g\gamma)^2$ are the same, contradicting the above remark.

Therefore the third condition is satisfied.

We come now to the fourth condition of a mixed Beauville quadruple. Let $g \in H_{[4]}$, $g \notin H_{[2]}$ be given, for instance $(1_H, 1_H, 1)$. Conjugation with $g$ interchanges then the first two components of an element $h \in H_{[4]}$. Our hypothesis 4 implies the result. □

As an application we find the following examples

**Theorem 28.** *Let $p$ be a prime with $p \equiv 3 \mod 4$ and $p \equiv 1 \mod 5$ and consider the group $H := \mathbf{SL}(2, \mathcal{F}_p)$. Then $H_{[4]}$ admits a mixed Beauville structure $u$ such that $S(u)$ is not biholomorphic to $\overline{S(u)}$.*

*Remark 11.* Note that the smallest prime satifying the above congruences is $p = 11$ and we get that $G$ has order equal to 6969600.

*Question 6.* : which is the minimal order of a group admitting a mixed Beauville structure?

## 5 Lefschetz pencils and braid monodromies

### 5.1 Braids and the mapping class group

The elegant definition by E. Artin of the *braid group* (cf. [Art26], [Art65]) supplies a powerful tool, even if difficult to handle, for the study of the differential topology of algebraic varieties, in particular of algebraic surfaces.

*Remark 12.* We observe that the subsets $\{w_1, \ldots, w_n\} \subset \mathbb{C}$ of $n$ distinct points in $\mathbb{C}$ are in one to one correspondence with monic polynomials $P(z) \in \mathbb{C}[z]$ of degree $n$ with non vanishing discriminant $\delta(P)$.

**Definition 16.** *The group*

$$\mathcal{B}_n := \pi_1(\mathbb{C}[z]_n \backslash \{P | \delta(P) = 0\}),$$

*i.e., the fundamental group of the space of polynomials of degree $n$ having $n$ distinct roots, is called* Artin's braid group.

Usually, one takes as base point the polynomial $P(Z) = (\prod_{i=1}^n (z-i)) \in \mathbb{C}[z]_n$ (or the set $\{1, \ldots, n\}$).

To a closed (continuous) path $\alpha : [0,1] \to (\mathbb{C}[z]_n \backslash \{P | \delta(P) = 0\})$ one associates the subset $\{(z,t) \in \mathbb{C} \times \mathbb{R} \mid \alpha(t)(z) := \alpha_t(z) = 0\}$ of $\mathbb{R}^3$.

Figure 2 below shows two realizations of the same braid.

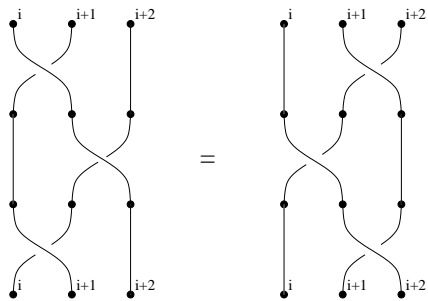

**Fig. 2.** Relation $aba = bab$ on braids

*Remark 13.* Obviously there is a lifting of $\alpha$ to $\mathbb{C}^n$, the space of $n$-tuples of roots of polynomials of degree $n$ and there are (continuous) functions $w_i(t)$ such that $w_i(0) = i$ and $\alpha_t(z) = \prod_{i=1}^n (z - w_i(t))$.

Then to each braid is associated a naturally defined permutation $\tau \in \mathfrak{S}_n$ given by $\tau(i) := w_i(1)$.

A very powerful generalization of Artin's braid group was given by M. Dehn (cf. [Deh38], we refer also to the book [Bir74]).

**Definition 17.** *Let $M$ be a differentiable manifold, then the* mapping class group *(or* Dehn group*) of $M$ is the group*

$$Map(M) := \pi_0(Diff(M)) = (Diff(M)/Diff^0(M)),$$

*where $Diff^0(M)$ is the subgroup of diffeomorphisms of $M$ isotopic to the identity.*

*Remark 14.* If $M$ is oriented then we often tacitly take $Diff^+(M)$, the group of orientation preserving diffeomorphisms of $M$ instead of $Diff(M)$, in the definition of the mapping class group. But it is more accurate to distinguish in this case $Map^+(M)$ from $Map(M)$.

If $M$ is a curve of genus $g$, then its mapping class group will be denoted by $Map_g$.

The relation between the above two definitions is the following:

**Theorem 29.** *The braid group $\mathcal{B}_n$ is isomorphic to the group*

$$\pi_0(Map^\infty(\mathbb{C}\setminus\{1,\ldots n\})),$$

*where $Map^\infty(\mathbb{C}\setminus\{1,\ldots n\})$ is the group of diffeomorphisms which are the identity outside the circle with center $0$ and radius $2n$.*

Therefore Artin's standard generators $\sigma_i$ of $\mathcal{B}_n$ ($i = 1,\ldots n-1$) can be represented by so-called half-twists.

**Definition 18.** *The half-twist $\sigma_j$ is the diffeomorphism of $\mathbb{C}\setminus\{1,\ldots n\}$ isotopic to the homeomorphism given by:*
*- rotation of $180$ degrees on the circle with center $j + \frac{1}{2}$ and radius $\frac{1}{2}$,*
*- on a circle with the same center and radius $\frac{2+t}{4}$ the map $\sigma_j$ is the identity if $t \geq 1$ and rotation of $180(1-t)$ degrees, if $t \leq 1$.*

Now, it is obvious that $\mathcal{B}_n$ acts on the free group $\pi_1(\mathbb{C}\setminus\{1,\ldots n\})$, which has a geometric basis (we take as base point the complex number $p := -2ni$) $\gamma_1,\ldots \gamma_n$ as explained in figure 3.

This action is called the *Hurwitz action of the braid group* and has the following algebraic description

- $\sigma_i(\gamma_i) = \gamma_{i+1}$
- $\sigma_i(\gamma_i\gamma_{i+1}) = \gamma_i\gamma_{i+1}$, whence $\sigma_i(\gamma_{i+1}) = \gamma_{i+1}^{-1}\gamma_i\gamma_{i+1}$
- $\sigma_i(\gamma_j) = \gamma_j$ for $j \neq i, i+1$.

Observe that the product $\gamma_1\gamma_2\ldots\gamma_n$ is left invariant under this action.

**Definition 19.** *We consider a group $G$ and its cartesian product $G^n$. The map associating to each $(g_1, g_2, \ldots, g_n)$ the product $g := g_1g_2\ldots, g_n \in G$ gives a partition of $G^n$, whose subsets are called* factorizations *of an element $g \in G$.*

*$\mathcal{B}_n$ acts on $G^n$ leaving invariant the partitions, and its orbits are called Hurwitz equivalence classes of factorizations.*

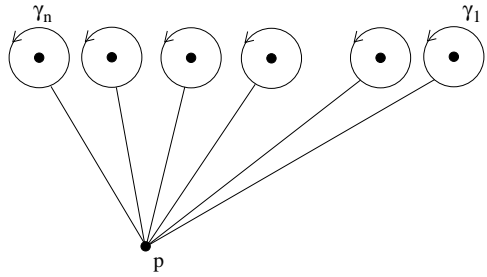

**Fig. 3.** A geometric base of $\pi_1(\mathbb{C} - \{1, \dots n\})$

**Definition 20.** *(cf. figure 4 below)*

*Let $C$ be a compact Riemann surface. Then a positive Dehn twist $T_\alpha$ with respect to a simple closed curve $\alpha$ on $C$ is an isotopy class of a diffeomorphism $h$ of $C$ which is equal to the identity outside an annular neighbourhood of $\alpha$, while inside the annulus $h$ rotates one boundary of the annulus by 360 degrees to the right and damps the rotation down to the identity at the other boundary.*

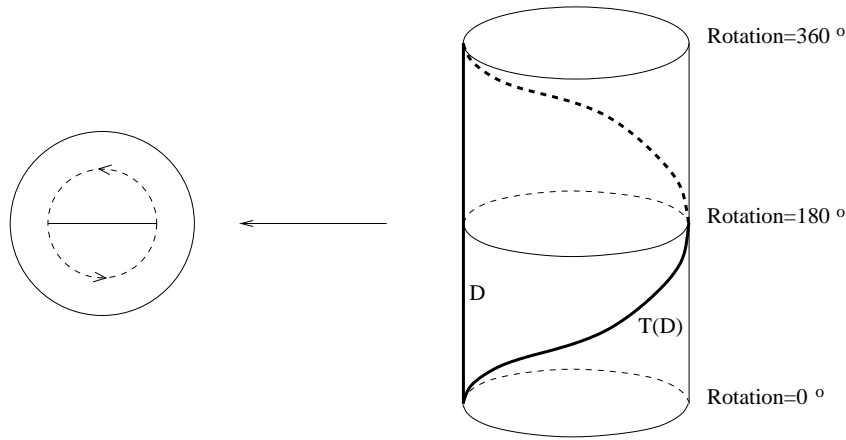

**Fig. 4.** At the left, a half twist; at the right: its lift: the Dehn twist $T$ and its action on the segment $D$

If one considers a hyperelliptic Riemann surface given as a branched cover of $\mathbb{P}^1_{\mathbb{C}} \setminus \{1, \dots n\}$ one sees that the Artin half twist $\sigma_j$ lifts to the Dehn twists on the loop which is the inverse image of the segment $[j, j+1]$.

Dehn's fundamental result is the following

**Theorem 30.** *The mapping class group $Map_g$ is generated by Dehn twists.*

Explicit presentations of $Map_g$ have been given by Hatcher and Thurston ([HT80]), which have been improved by Wajnryb ([Waj83]) who obtained a simpler presentation of the mapping class group (cf. also [Waj99]).

### 5.2 Lefschetz fibrations

The method introduced by Lefschetz for the study of the topology of algebraic varieties is the topological analogue of the method of hyperplane sections and projections of the classical italian algebraic geometers. It was classically used to describe the homotopy and homology groups of algebraic varieties. In the 70's Moisezon and Kas realized, after the work of Smale, that Lefschetz fibrations could be used to investigate the differential topology of algebraic varieties, especially of algebraic surfaces.

For instance, it is an extremely difficult problem to decide whether two algebraic surfaces which are not deformation equivalent are in fact diffeomorphic,even in the case where they are simply connected.

Here, the theory of Lefschetz fibrations offers a method to prove that two surfaces are diffeomorphic ([Kas80]).

**Definition 21.** *Let $M$ be a compact differentiable (or even symplectic) manifold of dimension* 4

*A* Lefschetz fibration *is a differentiable map $f : M \to \mathbb{P}^1_{\mathbb{C}}$ which*

*a) is of maximal rank except for a finite number of critical points $p_1, \ldots p_m$ which have distinct critical values $b_1, \ldots b_m \in \mathbb{P}^1_{\mathbb{C}}$,*

*b) has the property that around $p_i$ there are complex coordinates $(x, y) \in \mathbb{C}^2$ such that locally $f = x^2 - y^2 + const.$ (in the symplectic case, in the given coordinates the symplectic form $\omega$ of $M$ has to correspond to the natural symplectic structure on $\mathbb{C}^2$).*

*Remark 15.* 1) A similar definition can be given if $M$ is a manifold with boundary, replacing $\mathbb{P}^1_{\mathbb{C}}$ by a disc $D \subset \mathbb{C}$.

2) An important theorem of Donaldson ([Don99]) asserts that for symplectic manifolds there exists (as for the case of projective manifolds) a Lefschetz pencil, i.e., a Lefschetz fibration $f : M' \to \mathbb{P}^1_{\mathbb{C}}$ on a symplectic blow up $M'$ of $M$ (cf.[MS98]).

3) A Lefschetz fibration with fibres genus $g$ curves and with critical values $b_1, \ldots b_m \in \mathbb{P}^1_{\mathbb{C}}$, once a geometric basis $\gamma_1, \gamma_2, \ldots, \gamma_m$ of $\pi_1(\mathbb{P}^1_{\mathbb{C}} \setminus \{b_1, \ldots, b_m\})$ is chosen, determines a factorization of the identity in the mapping class group $Map_g$

$$\tau_1 \circ \tau_2 \circ \cdots \circ \tau_m = Id$$

as a product of Dehn twists.

We are now ready to state the theorem of Kas (cf. [Kas80]).

**Theorem 31.** *Two Lefschetz fibrations $(M, f)$, $(M', f')$ are equivalent (i.e., there are two diffeomorphisms $u : M \to M', v : \mathbb{P}^1 \to \mathbb{P}^1$ such that $f' \circ u = v \circ f$) if and only if the two corresponding factorizations of the identity in the mapping class group are equivalent (under the equivalence relation generated by Hurwitz equivalence and by simultaneous conjugation).*

*Remark 16.* 1) A similar result holds for Lefschetz fibrations over the disc and we get a factorization
$$\tau_1 \circ \tau_2 \circ \cdots \circ \tau_m = \phi$$
of the monodromy $\phi$ of the fibration over the boundary of the disc $D$.

2) The fibration admits a symplectic structure if and only if each Dehn twist in the factorization is positively oriented.

Assume that we are given two Lefschetz fibrations over $\mathbb{P}^1_{\mathbb{C}}$: then we can define the fiber sum of these two fibrations, which depends on a diffeomorphism chosen between two respective smooth fibers (cf. [GS99]).

This operation translates (in view of the above quoted theorem of Kas) into the following definition of "conjugated composition" of factorization:

**Definition 22.** *Let $\tau_1 \circ \tau_2 \circ \cdots \circ \tau_m = \phi$ and $\tau'_1 \circ \tau'_2 \circ \cdots \circ \tau'_r = \phi'$ be two factorizations: then their by $\psi$ conjugated composition is the factorization*
$$\tau_1 \circ \tau_2 \circ \ldots \tau_m \circ (\tau'_1)_\psi \circ (\tau'_2)_\psi \circ \cdots \circ (\tau'_r)_\psi = \phi(\phi')_\psi.$$

*Remark 17.* 1) If $\psi$ and $\phi'$ commute, we obtain a factorization of $\phi\phi'$.

2) A particular case is $\phi, \phi' = id$ and it corresponds to Lefschetz fibrations over $\mathbb{P}^1$.

### 5.3 Braid monodromy and Chisini' problem

Let $B \subset \mathbb{P}^2_{\mathbb{C}}$ be a plane algebraic curve of degree $d$, and let $P$ be a generic point not on $B$. Then the pencil of lines $L_t$ passing through $P$ determines a one parameter family of $d$-uples of points of $\mathbb{C} \cong L_t\backslash\{P\}$, i.e., $L_t \cap B$. Therefore one gets a factorization of $(\Delta^2)^d$ in the braid group $\mathcal{B}_n$, where $(\Delta^2) = (\sigma_{d-1}\sigma_{d-2}\ldots\sigma_1)^d$ is the generator of the center of the braid group. The equivalence class of the factorization does not depend on the point $P$ (if it is chosen generic) and does not depend on $B$, if $B$ varies in an equisingular family of curves.

Chisini was mainly interested in the case of *cuspidal* curves (cf. e.g. [Chi44], [Chi55]), mainly because these are the branch curves of a generic projection $f : S \to \mathbb{P}^2_{\mathbb{C}}$, for any smooth projective surface $S \subset \mathbb{P}^r$.

More precisely, a *generic projection* $f : S \to \mathbb{P}^2_{\mathbb{C}}$ is a covering whose branch curve has only nodes and cusps as singularities, and moreover is such that the local monodromy around a smooth point of the branch curve is a transposition.

Maps with those properties are called *generic coverings*: for these the local monodromies are only $\mathbb{Z}/2 = \mathfrak{S}_2$ (at the smooth points of the branch curve $B$), $\mathfrak{S}_3$ at the cusps, and $\mathbb{Z}/2 \times \mathbb{Z}/2$ at the nodes.

In such a case we have a *cuspidal* factorization, i.e. all factors are powers of a half twist, with respective exponent $1, 2, 3$.

Chisini posed the following

*Conjecture 3.* (Chisini's conjecture.)

Given two generic coverings $f : S \to \mathbb{P}^2_{\mathbb{C}}$, $f' : S' \to \mathbb{P}^2_{\mathbb{C}}$, one of them of degree $d \geq 5$, assume that they have the same branch curve $B$. Is it then true that $f$ and $f'$ are equivalent?

Observe that the condition on the degree is necessary, since counterexamples with $d = 4$ are furnished by the dual curve of a smooth plane cubic (as already known to Chisini, who gave a counterexample with $d = 4, d' = 3$, while counterexamples with $d = d' = 4$ were given in [Cat86b]).

The conjecture has been proven under the hypothesis that the degree of each covering is at least 12, essentially by Kulikov (cf. [Kul99]). In fact, Kulikov proved the result under a more complicated assumption and shortly later Nemirovski [Nem01] noticed, just by using the Miyaoka-Yau inequality, that Kulikov's assumption was implied by the simple assumption $d \geq 12$. Later on generalizations of this result were obtained for singular (normal) surfaces [Kul03] or for curves with more complicated singularities [MP02].

A negative answer instead has the following problem of Chisini (due to work of B. Moishezon (cf. [Moi94]).

**Chisini' s problem:** (cf.[Chi55]).

Given a cuspidal factorization, which is regenerable to the factorization of a smooth plane curve, is there a cuspidal curve which induces the given factorization?

*Regenerable* means that there is a factorization (in the equivalence class) such that, after replacing each factor $\sigma^i$ ($i = 2, 3$) by the $i$ corresponding factors (e.g. , $\sigma^3$ is replaced by $\sigma \circ \sigma \circ \sigma$) one obtains the factorization belonging to a non singular plane curve.

*Remark 18.* 1) Moishezon proves that there exist infinitely many non equivalent cuspidal factorizations observing that $\pi_1(\mathbb{P}^2_{\mathbb{C}} \setminus B)$ is an invariant defined in terms of the factorization alone. On the other hand, the family of cuspidal curves of a fixed degree form an algebraic set, hence has a finite number of connected components. These two statements together give a negative answer to the above cited problem of Chisini.

The examples of Moishezon have been recently reinterpreted in [ADK03], with a simpler treatment, in terms of symplectic surgeries.

2) In fact, as conjectured by Moishezon, a cuspidal factorization together with a generic monodromy with values in $\mathfrak{S}_n$ induces a covering $M \to \mathbb{P}^2_{\mathbb{C}}$, where $M$ is a symplectic fourmanifold.

Extending Donaldson's techniques (for proving the existence of symplectic Lefschetz fibrations) Auroux and Katzarkov ([AK00]) proved that each symplectic 4-manifold is in a natural way 'asymptotically' realized by such a generic covering. They propose to use an appropriate quotient of $\pi_1(\mathbb{P}^2_{\mathbb{C}} \backslash B)$ in order to produce invariants of symplectic structures, using the methods introduced by Moishezon and Teicher in a series of technically difficult papers ( see e.g. [MT92]).

It seems however that, up to now, these groups $\pi_1(\mathbb{P}^2_{\mathbb{C}} \backslash B)$ allow only to detect homology invariants of the projected fourmanifold ([ADKY04]).

3) Suppose we have a surface $S$ of general type and a pluricanonical embedding. Then a generic projection to $\mathbb{P}^3_{\mathbb{C}}$ gives a surface with a double curve $\Gamma'$. Now, project further to $\mathbb{P}^2_C$ and we do not only get the branch curve $B$, but also a curve $\Gamma$, image of $\Gamma'$.

Even if Chisini's conjecture tells us that from the holomorphic point of view $B$ determines the surface $S$ and therefore the curve $\Gamma$, it does not follow that the fundamental group $\pi_1(\mathbb{P}^2_{\mathbb{C}} \backslash B)$ determines the group $\pi_1(\mathbb{P}^2_{\mathbb{C}} \backslash (B \cup \Gamma))$.

It would be interesting to calculate this second fundamental group, even in special cases.

## 6 DEF, DIFF and other equivalence relations

As we said, one of the fundamental problems in the theory of complex algebraic surfaces is to understand the moduli spaces of surfaces of general type, and in particular their connected components, which parametrize the deformation equivalence classes of minimal surfaces of general type.

**Definition 23.** *Two minimal surfaces $S$ and $S'$ are said to be* def *- equivalent (we also write: $S \sim_{def} S'$) if and only if they are elements of the same connected component of the moduli space.*

By the classical theorem of Ehresmann, two def - equivalent algebraic surfaces are (orientedly) diffeomorphic.

In the late eighties Friedman and Morgan (cf. [FM88]) conjectured that two algebraic surfaces are diffeomorphic if and only if they are def - equivalent. We will abbreviate this conjecture in the following by the acronym **def = diff**.

The second author would like to point out here that he had made the opposite conjecture in the early eithties (cf. [Kat83]).

Donaldson's breaktrough results had made clear that diffeomorphism and homeomorphism differ drastically for algebraic surfaces (cf. [Don83]) and the success of gauge theory led Frieman and Morgan to "speculate" that the diffeomorphism type of algebraic surfaces determines the deformation class. After the first counterexamples of M. Manetti (cf. [Man01]) appeared, there were further counterexamples given by Catanese, Kharlamov-Kulikov, Catanese-Wajnryb, Bauer-Catanese-Grunewald (cf. [Cat03], [KK02],[BCG05a], [CW04]).

In the cited papers by Catanese, Kharlamov-Kulikov, Bauer-Catanese-Grunewald, the counterexamples are given by pairs of surfaces, where one is the complex conjugate of the other.

One could say that somehow these counterexamples are 'cheap', and somehow in the air (cf. the definition of strong rigidity). The second author was very recently informed by R. Friedman that also he and Morgan were aware of such 'complex conjugate' counterexamples, but for the case of elliptic surfaces.

Since the beautiful examples of Manetti yield non simply connected surfaces, it made sense to weaken the conjecture **def = diff** in the following way.

*Question 7.* Is the speculation **def = diff** true if one requires the diffeomorphism $\phi : S \to S'$ to send the first Chern class $c_1(K_S) \in H^2(S, \mathbb{Z})$ in $c_1(K_{S'})$ and moreover one requires the surfaces to be simply connected?

But even this weaker question turned out to have a negative answer, as it was shown by the second author and Wajnryb ([CW04]).

*Remark 19.* If two surfaces are def-equivalent, then there exists a diffeomrophism sending the canonical class $c_1(K_S) \in H^2(S, \mathbb{Z})$ in the canonical class $c_1(K_{S'})$. On the other hand, by the result of Seiberg - Witten theory we know that a diffeomorphism sends the canonical class of a minimal surface $S$ to $\pm c_1(K_{S'})$. Therefore, if one gives at least three surfaces, which are pairwise diffeomorphic, one finds at least two surfaces with the property that there exists a diffeomorphism between them sending the canonical class of one to the canonical class of the other.

**Theorem 32.** *([CW04])*

*For each natural number h there are simply connected surfaces $S_1, \ldots, S_h$ which are pairwise diffeomorphic, but are such that two of them are never def - equivalent.*

The above surfaces $S_1, \ldots, S_h$ belong to the class of the so-called $(a, b, c)$-surfaces, obtained as minimal compactification of some affine surface described by the following two equations:

$$z^2 = f(x, y),$$
$$w^2 = g(x, y),$$

where $f$ and $g$ are suitable polynomials of respective bidegrees $(2a, 2b)$, $(2c, 2b)$.

They can be compactified simply by bihomogenizing the polynomials $f, g$, thus obtaining Galois covers of $\mathbb{P}^1 \times \mathbb{P}^1$ with Galois group $\mathbb{Z}/2\mathbb{Z}$.

We remark that the above compactification is smooth if the two curves $\{f = 0\}$ and $\{g = 0\}$ in $\mathbb{P}^1 \times \mathbb{P}^1$ are smooth and intersect transversally.

We say that these surfaces are bidouble (i.e., Galois covers with Galois group $(\mathbb{Z}/2\mathbb{Z})^2$) covers of $\mathbb{P}^1 \times \mathbb{P}^1$ of type $(2a, 2b)$, $(2c, 2b)$ (cf. [Cat84], [Cat99]).

The above theorem is implied by the two following results:

**Theorem 33.** *Let $a, b, c, k$ be positive even numbers such that*
  *1) $a, b, c - k \geq 4$;*
  *2) $a \geq 2c + 1$;*
  *3) $b \geq c + 2$.*
  *Furthermore, let $S$ be an $(a, b, c)$ - surface and $S'$ be an $(a + k, b, c - k)$-surface. Then $S$ is not def - equivalent to $S'$.*

**Theorem 34.** *Let $S$ be an $(a, b, c)$ - surface and $S'$ be an $(a + 1, b, c - 1)$-surface. Moreover, assume that $a, b, c - 1 \geq 2$. Then $S$ and $S'$ are diffeomorphic.*

*Remark 20.* Observe that the surfaces in question above are simply connected (cf. [Cat84], prop. 2.7.).

The proof of the two theorems above are completely different in nature. The first theorem uses techniques which have been developed in a series of papers by the second author and by Manetti ([Cat84], [Cat87a], [Cat86a], [Man94], [Man97]). They use essentially the local deformation theory a la Kuranishi, normal degenerations of smooth surfaces and a study of quotient singularities of rational double points and of their smoothings.

One very elementary, but extremely important ingredient in the proof of the first theorem is the notion of *natural deformations of a bidouble cover* (introduced in [Cat84], p.494), which are parametrized by a quadruple of polynomials $(f, g, \phi, \psi)$ and given by the two equations

$$z^2 = f(x, y) + w\phi(x, y),$$

$$w^2 = g(x, y) + z\psi(x, y),$$

where $f$ and $g$ are polynomials of respective bidegrees $(2a, 2b)$, $(2c, 2b)$ as before and $\phi$ and $\psi$ have respective bidegrees $(2a - c, b)$, $(2c - a, b)$.

Under suitable hypotheses (rigid base, branch curve of sufficiently high degree), these natural deformations indeed give all small deformations. Moreover, since $a \geq 2c + 1$, it follows that $\psi \equiv 0$, therefore every small deformation preserves the structure of an iterated double cover.

The final point is to show that this structure also passes in a suitable way to the limit, so that we do not only have an open, but also a closed subset of the moduli space.

We will now comment on the newer part, the proof of theorem 34.

The key ideas are here the following:

1) Both surfaces $S$ and $S'$ admit a holomorphic map to $\mathbb{P}^1_{\mathbb{C}}$ given by the composition of the bidouble cover with the projection to the first coordinate $x$, and a small perturbation of this map realizes them as symplectic Lefschetz fibrations (cf. [Don99], [GS99]).

2) The respective fibrations are, by the accurate choice of the bidegrees of the curves, and especially because of the fact that the second degree is equal to

46    Ingrid C. Bauer, Fabrizio Catanese, and Roberto Pignatelli$2b$ in both cases (thereby allowing, locally on the base, to 'rotate' one branch curve to the other) fiber sums of the same pair of Lefschetz fibrations over the complex disc (the global effect of this local rotation is that the first curve $\{f = 0\}$ loses a bidegree $(2,0)$, while the second $\{g = 0\}$ gains a bidegree $(2,0)$))

3) Once the first fiber sum is presented as composition of two factorizations and the second as the same composition of factorizations, just conjugated by the 'rotation' $\Psi$, in order to prove that the two fiber sums are equivalent, it suffices, (thanks to Auroux's lemma, [Aur02]) to show that the diffeomorphism $\Psi$ is in the subgroup of the mapping class group generated by the Dehn twists which appear in the first factorization.

4) Figure 5 below shows the fibre $C$ of the fibration in the case $2b = 6$: it is a bidouble cover of $\mathbb{P}^1$, which we can assume to be given by the equations $z^2 = F(y)$, $w^2 = F(-y)$, where the roots of $F$ are the integers $1, \ldots, 2b$.

Moreover, one sees that the monodromy of the fibration at the boundary of the disc is trivial, and the map $\Psi$ is the diffeomorphism of order 2 given by $y \mapsto -y$, $z \mapsto w$, $w \mapsto z$, which in our figure is given as a rotation of 180 degrees around an axis inclined in direction north-east.

The figure shows a dihedral symmetry, where the automorphism of order 4 is given by $y \mapsto -y$, $z \mapsto -w$, $w \mapsto z$.

Moreover, between the Dehn twists which appear in the factorization there are those which correspond to the inverse images of the segments between two consecutive integers (cf. figure 5). These circles can be organized on the curve $C$ in six chains (not disjoint) and finally we have reduced ourselves to show that the isotopy class of $\Psi$ is the same as the product of the six Coxeter elements associated to such chains.

We recall that the *Coxeter elements associated to a chain* are products of the type
$$\Delta = (T_{\alpha_1})(T_{\alpha_2}T_{\alpha_1}) \ldots (T_{\alpha_n}T_{\alpha_{n-1}} \ldots T_{\alpha_1})$$
of Dehn twists associated to the curves of the chain.

In order to finally prove that such product (let us call it $\Psi'$) of Coxeter elements and $\Psi$ are isotopic, one observes that if one removes the above cited chains of circles from the curve $C$, one obtains 4 connected components which are diffeomorphic to circles. By a result of Epstein it is then sufficient to verify that $\Psi$ and $\Psi'$ send each such curve to a pair of isotopic curves: this last step needs a list of lengthy (though easy) verifications, for which it is necessary to have explicit drawings.

For details we refer to the original paper [CW04].

It was observed by the second author (cf. [Cat02]) that a surface of general type has a canonical symplectic structure. In fact, he proves the following

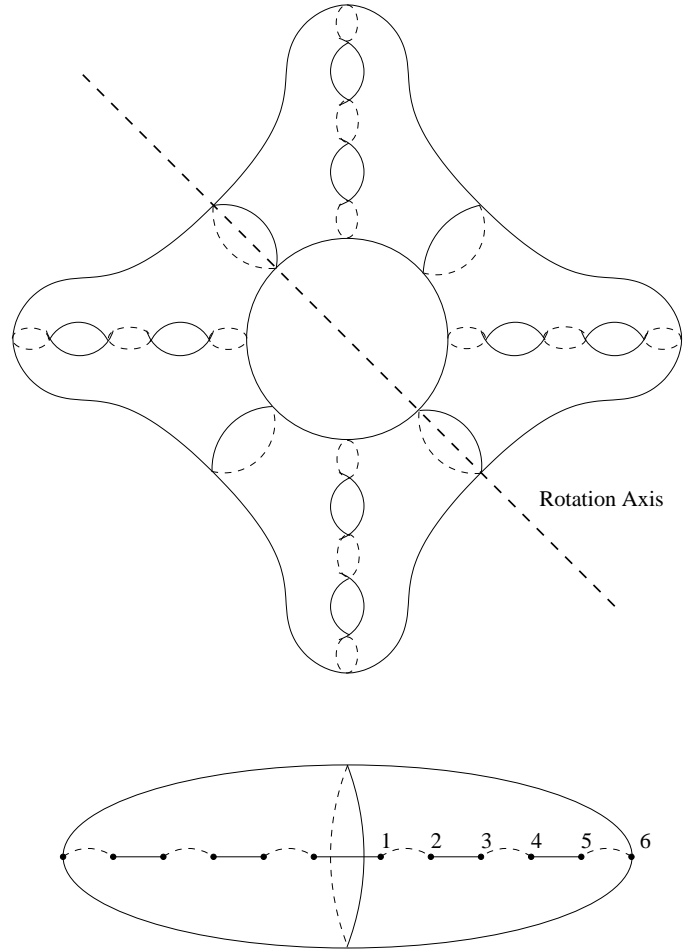

**Fig. 5.** The curve $C$ with a dihedral symmetry

**Theorem 35.** *A minimal surface of general type $S$ has a canonical symplectic structure, unique up to symplectomorphism, such that the class of the symplectic form is the class of the canonical sheaf $\mathcal{O}_S(K_S)$.*

We give the proof under the assumption that $K_S$ be ample, since the basic idea becomes clear in this simpler case. For the general case we refer to the original article.

*Proof.* Let $m$ be such that $mK_S$ is very ample (by Bombieri's result it suffices any $m \geq 4$), i.e., the pluricanonical map $\phi_m := \phi_{mK_S} : S \to \mathbb{P}^{P_m-1}$, where $P_m := h^0(S, \mathcal{O}_S(mK_S))$ is the $m$-th plurigenus of $S$, is an embedding.

We define $\omega_m$ on $S$ as follows:

$$\omega_m := \frac{1}{m}\phi_m^*(\frac{1}{2\pi i}\partial\overline{\partial}log|z|^2),$$

i.e., we divide by $m$ the pull back of the Fubini - Study form, whence $\omega_m$ yields a symplectic form on $S$.

It remains to show that the symplectomorphism class of $(S, \omega_m)$ is indeed independent of $m$.

For this suppose that also $\phi_n$ gives an embedding of $S$: then the same holds for $mn$, whence it is sufficient to see that $(S, \omega_m)$ and $(S, \omega_{nm})$ are symplectomorphic. Observe that the pull back of the Fubini-Study form under the $n$-th Veronese map $v_n$ is $n$ times the Fubini-Study form and $v_n \circ \phi_m$ is a linear projection of $\phi_{mn}$. Then by Moser's theorem we are done. □

Therefore it seems natural to ask the following

*Question 8.* Are the diffeomorphic $(a, b, c)$-surfaces of theorem 34, endowed with their canonical symplectic structure, indeed symplectomorphic?

*Remark 21.* 1) In [Cat02] the second author shows that Manetti's examples are indeed symplectomorphic.

2) A possible way of showing that the answer to the question above is yes (and therefore exhibiting symplectomorphic simply connected surfaces which are not def-equivalent) goes through the analysis of the braid monodromy of the branch curve of the "perturbed" (corresponding to the Lefschetz fibration) quadruple covering, and one would like to show that the involution $\iota$ on $\mathbb{P}^1$, $\iota(y) = -y$ can be written as the product of braids which show up in the factorization.

Anyhow, this approach turned out to be more difficult than the corresponding analysis which has been made in the mapping class group, because the braid monodromy contains very many 'tangency' factors which do not come from local contributions to the regeneration of the branch curve from the union of the curves $f = 0, g = 0$ counted twice.

In the rest of the paragraph we will discuss another equivalence relation the so-called **Q.E.D.** equivalence relation, which was introduced by the second author (cf. [Cat05]), and which seems worthwhile to examine for surfaces of general type.

Observe that for any number $g \geq 2$ there is a smooth curve of genus $g$, which is an étale covering of a curve of genus 2. Therefore all the smooth curves of Kodaira dimension 1 are equivalent by the equivalence relation generated by deformation and by étale maps. (Obviously, also all curves of Kodaira dimension 0, resp. $-\infty$ are equivalent by this equivalence relation).

*Remark 22.* More remarkable is what happens for algebraic surfaces of Kodaira dimension 0. Enriques surfaces admit an étale double cover which is a K3-surface, hyperelliptic surfaces have an étale cover which is a torus (in fact, the product of two elliptic curves).

Therefore, in order to have some analogue to the curve case, one should "link" K3-surfaces and tori by étale maps and deformations. Obviously, this is not possible, since tori are $K(\pi,1)$'s and K3-surfaces are simply connected.

But the solution is simple: divide the torus by the involution $x \mapsto -x$, and obtain the (singular!) Kummer surface. A smoothing of this Kummer surface gives a K3-surface. The price we have to pay for going from curves to surfaces is that we have to allow morphisms which are not necessarily étale, but only étale in codimension 1. Moreover, we have to allow mild singularities: ordinary double points in this case, canonical singularities in a more general setting.

This remark justifies the following

**Definition 24.** *We consider for complete algebraic varieties with canonical singularities defined over a fixed algebraically closed field the equivalence relation generated by*

*1) birational maps;*

*2) flat proper algebraic deformations $\pi : \mathcal{X} \to B$, with base $B$ a connected algebraic variety, and all fibres having canonical singularities;*

*3) quasi étale morphisms $f : X \to Y$, i.e., surjective morphisms which are étale in codimension 1 on $X$ (i.e., there is $\mathcal{Z} \subset X$ of codimension $\geq 2$ such that $f|(X - Z)$ is étale).*

*We will call this equivalence relation* a.q.e.d. - relation*, which means algebraic quasiétale - deformation relation and it will be denoted by $X \sim_{a.q.e.d.} X'$.*

It is rather clear that a completely analogous equivalence relation (called then $\mathbb{C}$ - *q.e.d.* - *relation*) can be defined also in the setting of compact complex spaces with canonical singularities. We refer to [Cat05] for more details.

*Remark 23.* Trivially, the dimension of a variety is a q.e.d. invariant.

By Siu's recent result (cf. [Siu02]) also the Kodaira dimension is an invariant of a.q.e.d. - equivalence, if we restrict ourselves to projective varieties with canonical singularities (defined over the complex numbers).

For surfaces of special type, i.e., of Kodaira dimension $\leq 1$ the situation is as for curves.

**Theorem 36.** *Let $S$ and $S'$ be smooth complex algebraic surfaces of the same Kodaira dimension $\leq 1$. Then $S$ and $S'$ are a.q.e.d. - equivalent.*

The ingredients of the proof of the above theorem are the Enriques classification of surfaces, the detailed knowledge of the deformation types of elliptic surfaces and the orbifold fundamental group of a fibration.

The following question seems natural.

*Question 9.* Is it possible to determine the q.e.d. equivalence classes inside the class of varieties with fixed dimension $n$, and with Kodaira dimension $k$?

For curves and special algebraic surfaces over $\mathbb{C}$ there is only one a.q.e.d. class, but as shown in an appendix to [Cat05] by Fritz Grunewald, already for surfaces of general type the situation is completely different.

**Theorem 37.** *There are infinitely many q.e.d. - equivalence classes of algebraic surfaces of general type.*

The above surfaces are constructed from quaternion algebras (along general lines suggested by Shimura and explicitly described by Kuga and Shavel, cf. [Sha78]) are rigid, but the q.e.d. - equivalence class contains countably many distinct birational classes.

The main points of the construction are the following:

1) the surfaces are quotients $S = \mathbb{H} \times \mathbb{H}/\Gamma$ of the two dimensional polydisk $\mathbb{H} \times \mathbb{H}$ via the free action of a discrete group $\Gamma$ constructed from a quaternion algebra $\mathcal{A}$ over a totally real quadratic field $k$

2) since $S$ is rigid, it suffices to show that if $\Gamma'$ is commensurable with $\Gamma$, then also $\Gamma'$ acts freely on $\mathbb{H} \times \mathbb{H}$ .

3) One sees by general theorems that $\Gamma'$ has as $\mathbb{Q}$-linear span the same quaternion algebra $\mathcal{A}$ as $\Gamma$.

4) If $\Gamma'$ does not act freely, taking the tangent representation at a fixed point, we see by 3) that $\mathcal{A}$ contains a cyclotomic extension whose degree divides 4.

5) Using Hasse's theorem, one chooses $\mathcal{A}$ such that the set of primes where it ramifies contains, one for each possible intermediate field $K'$ between the quadratic field $k$ of $\mathcal{A}$ and one of the finitely many possible cyclotomic extensions above, a prime $P$ such that $K' \otimes k_P$ is not an integral domain: this however contradicts 4) hence shows the desired assertion.

It remains open whether there are for instance varieties which are isolated in their q.e.d.-equivalence class (up to birational equivalence, of course).

An interesting question is to determine, for surfaces of general type, the non standard a.q.e.d. classes (standard means: equivalent to a product of curves)

$M_1$　　　　$M_2$

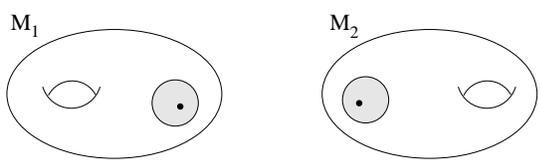

$M_1 \# M_2$

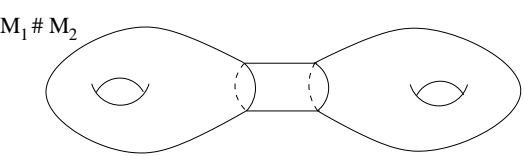

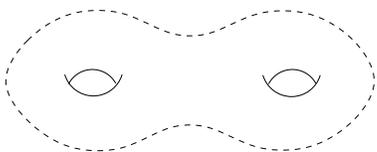

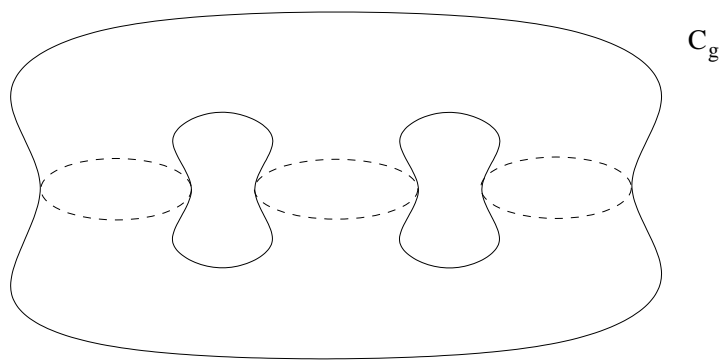

$C_g$

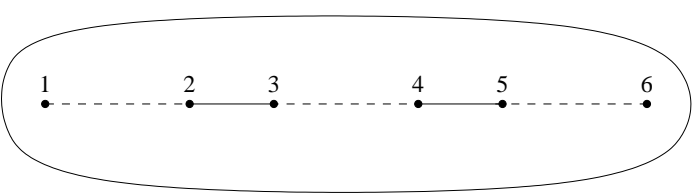

$P^1$